# MULTI–PEAK SOLUTIONS FOR A CLASS OF DEGENERATE ELLIPTIC EQUATIONS

ALESSANDRO GIACOMINI AND MARCO SQUASSINA

ABSTRACT. By means of a penalization argument due to del Pino and Felmer, we prove the existence of multi–spike solutions for a class of quasilinear elliptic equations under natural growth conditions. Compared with the semilinear case some difficulties arise, mainly concerning the properties of the limit equation. The study of concentration of the solutions requires a somewhat involved analysis in which a Pucci–Serrin type identity plays an important role.

## 1. INTRODUCTION AND MAIN RESULTS

Let $\Omega$ be any smooth domain of $\mathbb{R}^N$ with $N \geqslant 3$. Starting from the celebrated paper by Floer and Weinstein [15], much interest has been directed in recent years to the singularly perturbed semilinear elliptic problem

$$\begin{cases} -\varepsilon^2 \Delta u + V(x)u = u^{q-1} & \text{in } \Omega \\ u > 0 & \text{in } \Omega \\ u = 0 & \text{on } \partial\Omega \end{cases}$$

where $2 < q < 2N/(N-2)$ and $V(x)$ is a positive function bounded from below away from zero. Typically, there exists a family of solutions $(u_\varepsilon)_{\varepsilon>0}$ which exhibits a spike–shaped profile around every possibly degenerate local minimum of $V(x)$ and decade elsewhere as $\varepsilon$ goes to zero (see e.g. [1, 11, 13, 22, 25, 32] for the single–peak case and [12, 23] for the multi–peak case). A natural question is now whether these concentration phenomena are a special feature of the semilinear case or we can expect a similar behavior to hold for more general elliptic equations having a variational structure. The results concerning the existence of one–peak solutions have been recently

*Date*: 28th October 2018.
1991 *Mathematics Subject Classification.* 35J40; 58E05.
The research of the second author was partially supported by the MIUR project "Variational and topological methods in the study of nonlinear phenomena" (COFIN 2001) and by Gruppo Nazionale per l'Analisi Matematica, la Probabilità e le loro Applicazioni (INdAM).

extended in [30] to the quasilinear elliptic equation

$$-\varepsilon^2 \sum_{i,j=1}^{N} D_j(a_{ij}(x,u)D_i u) + \frac{\varepsilon^2}{2} \sum_{i,j=1}^{N} D_s a_{ij}(x,u)D_i u D_j u + V(x)u = u^{q-1}.$$

In this paper we turn to a more delicate situation, namely the study of the multi–peak case, also for possibly degenerate operators.

Assume that $V : \mathbb{R}^N \to \mathbb{R}$ is a $C^1$ function and there exists a positive constant $\alpha$ such that

(1) $$V(x) \geqslant \alpha \quad \text{for every } x \in \mathbb{R}^N.$$

Moreover let $\Lambda_1, \ldots, \Lambda_k$ be $k$ disjoint compact subsets of $\Omega$ and $x_i \in \Lambda_i$ with

(2) $$V(x_i) = \min_{\Lambda_i} V < \min_{\partial \Lambda_i} V, \quad i = 1, \ldots, k.$$

Let us set for all $i = 1, \ldots, k$

(3) $$\mathscr{M}_i := \{x \in \Lambda_i : V(x) = V(x_i)\}.$$

Let $1 < p < N$, $p^* = \frac{Np}{N-p}$ and let $W_V(\Omega)$ be the weighted Banach space

$$W_V(\Omega) := \left\{ u \in W_0^{1,p}(\Omega) : \int_\Omega V(x)|u|^p < +\infty \right\}$$

endowed with the natural norm $\|u\|_{W_V}^p := \int_\Omega |Du|^p + \int_\Omega V(x)|u|^p$. For all $A, B \subseteq \mathbb{R}^N$, let us denote their distance by $\text{dist}(A, B)$.

The following is the first of our main results.

**Theorem 1.1.** *Assume that* (1) *and* (2) *hold and let* $1 < p \leqslant 2$, $p < q < p^*$. *Then there exists* $\varepsilon_0 > 0$ *such that, for every* $\varepsilon \in (0, \varepsilon_0)$, *there exist* $u_\varepsilon$ *in* $W_V(\Omega) \cap C_{\text{loc}}^{1,\beta}(\Omega)$ *and* $k$ *points* $x_{\varepsilon,i} \in \Lambda_i$ *satisfying the following properties:*

(a) $u_\varepsilon$ *is a weak solution of the problem*

(4) $$\begin{cases} -\varepsilon^p \Delta_p u + V(x)u^{p-1} = u^{q-1} & \text{in } \Omega \\ u > 0 & \text{in } \Omega \\ u = 0 & \text{on } \partial \Omega; \end{cases}$$

(b) *there exist* $\sigma, \sigma' \in ]0, +\infty[$ *such that for every* $i = 1, \ldots, k$ *we have*

$$u_\varepsilon(x_{\varepsilon,i}) = \sup_{\Lambda_i} u_\varepsilon, \quad \sigma < u_\varepsilon(x_{\varepsilon,i}) < \sigma', \quad \lim_{\varepsilon \to 0} \text{dist}(x_{\varepsilon,i}, \mathscr{M}_i) = 0$$

   *where* $\mathscr{M}_i$ *is as in* (3);

(c) *for every* $r < \min\{\text{dist}(\mathscr{M}_i, \mathscr{M}_j) : i \neq j\}$ *we have*

$$\lim_{\varepsilon \to 0} \|u_\varepsilon\|_{L^\infty(\Omega \setminus \bigcup_{i=1}^k B_r(x_{\varepsilon,i}))} = 0;$$

(d) *it results*

$$\lim_{\varepsilon \to 0} \|u_\varepsilon\|_{W_V} = 0.$$

*Moreover, if $k = 1$ the assertions hold for every $1 < p < N$.*

Actually, this result will follow by a more general achievement involving a larger class of quasilinear operators. Before stating it, we make a few assumptions.

Assume that $1 < p < N$, $f \in C^1(\mathbb{R}^+)$ and there exist $p < q < p^*$ and $p < \vartheta \leqslant q$ with

$$\text{(5)} \qquad \lim_{s \to 0^+} \frac{f(s)}{s^{p-1}} = 0, \quad \lim_{s \to +\infty} \frac{f(s)}{s^{q-1}} = 0,$$

$$\text{(6)} \qquad 0 < \vartheta F(s) \leqslant f(s)s \quad \text{for every } s \in \mathbb{R}^+,$$

where $F(s) = \int_0^s f(t)\, dt$ for every $s \in \mathbb{R}^+$.

The function $j(x, s, \xi) : \Omega \times \mathbb{R}^+ \times \mathbb{R}^N \to \mathbb{R}$ is continuous in $x$ and of class $C^1$ with respect to $s$ and $\xi$, the function $\{\xi \mapsto j(x, s, \xi)\}$ is strictly convex and $p$–homogeneous and there exist two positive constants $c_1, c_2$ with

$$\text{(7)} \qquad |j_s(x, s, \xi)| \leqslant c_1 |\xi|^p, \quad |j_\xi(x, s, \xi)| \leqslant c_2 |\xi|^{p-1}$$

for a.e. $x \in \Omega$ and every $s \in \mathbb{R}^+$, $\xi \in \mathbb{R}^N$ ($j_s$ and $j_\xi$ denote the derivatives of $j$ with respect of $s$ and $\xi$ respectively). Let $R, \nu > 0$ and $0 < \gamma < \vartheta - p$ with

$$\text{(8)} \qquad j(x, s, \xi) \geqslant \nu |\xi|^p,$$

$$\text{(9)} \qquad j_s(x, s, \xi)s \leqslant \gamma j(x, s, \xi)$$

a.e. in $\Omega$, for every $s \in \mathbb{R}^+$ and $\xi \in \mathbb{R}^N$, and

$$\text{(10)} \qquad j_s(x, s, \xi) \geqslant 0 \quad \text{for every } s \geqslant R$$

a.e. in $\Omega$ and for every $\xi \in \mathbb{R}^N$. For every fixed $\bar{x} \in \Omega$, the limiting equation

$$\text{(11)} \quad -\operatorname{div}(j_\xi(\bar{x}, u, Du)) + j_s(\bar{x}, u, Du) + V(\bar{x})u^{p-1} = f(u) \quad \text{in } \mathbb{R}^N$$

admits a unique positive solution (up to translations). Finally, we assume that

$$\text{(12)} \qquad j(x_i, s, \xi) = \min_{x \in \Lambda_i} j(x, s, \xi), \quad i = 1, \ldots, k$$

for every $s \in \mathbb{R}^+$ and $\xi \in \mathbb{R}^N$, where the $x_i$s are as in (2).

We point out that assumptions (1), (2), (5) and (6) are the same as in [11, 12]. Conditions (7)–(10) are natural assumption, already used, for instance, in [2, 3, 4, 5, 28, 29].

The following result is an extension of Theorem 1.1.

**Theorem 1.2.** *Assume that (1), (2), (5), (6), (7), (8), (9), (10), (11), (12) hold.*

*Then there exists $\varepsilon_0 > 0$ such that, for every $\varepsilon \in (0, \varepsilon_0)$, there exist $u_\varepsilon$ in $W_V(\Omega) \cap C^{1,\beta}_{\text{loc}}(\Omega)$ and $k$ points $x_{\varepsilon,i} \in \Lambda_i$ satisfying the following properties:*

(a) $u_\varepsilon$ is a weak solution of the problem

$(P_\varepsilon)$
$$\begin{cases} -\varepsilon^p \text{div}\,(j_\xi(x,u,Du)) + \varepsilon^p j_s(x,u,Du) + V(x)u^{p-1} = f(u) & \text{in } \Omega \\ u > 0 & \text{in } \Omega \\ u = 0 & \text{on } \partial\Omega\,; \end{cases}$$

(b) there exist $\sigma, \sigma' \in\,]0, +\infty[$ such that for every $i = 1, \ldots, k$ we have

$$u_\varepsilon(x_{\varepsilon,i}) = \sup_{\Lambda_i} u_\varepsilon, \quad \sigma < u_\varepsilon(x_{\varepsilon,i}) < \sigma', \quad \lim_{\varepsilon \to 0} \text{dist}(x_{\varepsilon,i}, \mathcal{M}_i) = 0$$

where $\mathcal{M}_i$ is as in (3);

(c) for every $r < \min\{\text{dist}(\mathcal{M}_i, \mathcal{M}_j) : i \neq j\}$ we have

$$\lim_{\varepsilon \to 0} \|u_\varepsilon\|_{L^\infty(\Omega \setminus \bigcup_{i=1}^k B_r(x_{\varepsilon,i}))} = 0\,;$$

(d) it results

$$\lim_{\varepsilon \to 0} \|u_\varepsilon\|_{W_V} = 0.$$

Notice that if $k = 1$ assumption (11) can be dropped: in fact following the arguments of [30] it is possible to prove that the previous result holds without any uniqueness assumption, which instead, as in the semilinear case, seems to be necessary for the case $k > 1$. This holds true for the $p$−Laplacian problem (4) and for more general situation we refer the reader to [27].

Various difficulties arise in comparison with the semilinear framework (see also Section 5 of [30]). In order to study the concentration properties of $u_\varepsilon$ inside the $\Lambda_i$s (see Section 4), inspired by the recent work of Jeanjean and Tanaka [17], we make a repeated use of a Pucci–Serrin type identity [10] which has turned out to be a very powerful tool (see Section 3). It has to be pointed out that, in our possibly degenerate setting, we cannot hope to have $C^2$ solutions, but at most $C^{1,\beta}$ solutions (see [14, 31]). Therefore, the classical Pucci–Serrin identity [24] is not applicable in our framework. On the other hand, it has been recently shown in [10] that, under minimal regularity assumptions, the identity holds for locally Lipschitz solutions, provided that the operator is strictly convex in the gradient, which, from our viewpoint, is a very natural requirement (see Theorem 3.1). Under uniqueness assumptions this identity has also turned out to be useful in characterizing the exact energy level of the solution of (11). More precisely, we prove that (11) admits a least energy solution having the Mountain–Pass energy level (see Theorem 3.2).

We stress that the functionals we will study, associated with $(P_\varepsilon)$, are not even locally Lipschitz continuous (unless $j_s = 0$) and that tools of nonsmooth critical point theory will be emploied (see [6, 9, 16, 18]).

The plan of the paper is as follows:

– in Section 2, following the approach of [12], we construct the penalized functional $E_\varepsilon$ and we prove that it satisfies the $(CPS)_c$ condition (cf. Definition 6.6);
– in Section 3 we prove two important consequences of a Pucci–Serrin type variational identity (cf. Theorem 3.2 and Lemma 3.3);
– in Section 4 we study the concentration of the solutions $u_\varepsilon$ (cf. Lemmas 4.1, 4.2, 4.6);
– in Section 5, finally, we end up the proofs of Theorems 1.1 and 1.2;
– in Section 6 we quote some tools of nonsmooth critical point theory.

Notations:

$\|\cdot\|_p$ and $\|\cdot\|_{1,p}$ are the standard norms of the spaces $L^p(\Omega)$ and $W_0^{1,p}(\Omega)$;
$\mathrm{dist}(x, E)$ is the distance of $x$ from a set $E \subset \mathbb{R}^N$;
$\mathbb{1}_E$ is the characteristic function of a set $E \subset \mathbb{R}^N$;
$B_\varrho(y)$ is the ball in $\mathbb{R}^N$ of center $y$ and radius $\varrho > 0$; we set $B_\varrho := B_\varrho(0)$.

## 2. Penalization and $(CPS)_c$ condition

In this section, following the approach of del Pino and Felmer [12], we define a suitable penalization of the functional $I_\varepsilon : W_V(\Omega) \to \mathbb{R}$ associated with the problem $(P_\varepsilon)$,

$$I_\varepsilon(u) := \varepsilon^p \int_\Omega j(x, u, Du) + \frac{1}{p} \int_\Omega V(x) |u|^p - \int_\Omega F(u).$$

By the growth condition on $j$, it is easily seen that $I_\varepsilon$ is a continuous functional. We refer the reader to Section 6 for more details on the variational formulation in this nonsmooth framework. Let $\alpha > 0$ be as in (1) and consider the positive constant

$$(13) \qquad \ell := \sup \left\{ s > 0 : \frac{f(t)}{t^{p-1}} \leqslant \frac{\alpha}{\kappa} \text{ for every } 0 \leqslant t \leqslant s \right\}$$

for some fixed $\kappa > \vartheta/(\vartheta - p)$. We define the function $\widetilde{f} : \mathbb{R}^+ \to \mathbb{R}$ by setting

$$\widetilde{f}(s) := \begin{cases} \frac{\alpha}{\kappa} s^{p-1} & \text{if } s > \ell \\ f(s) & \text{if } 0 \leqslant s \leqslant \ell \end{cases}$$

and the map $g : \Omega \times \mathbb{R}^+ \to \mathbb{R}$ as

$$g(x, s) := \mathbb{1}_\Lambda(x) f(s) + (1 - \mathbb{1}_\Lambda(x)) \widetilde{f}(s), \quad \Lambda = \bigcup_{i=1}^k \Lambda_i$$

for a.e. $x \in \Omega$ and every $s \in \mathbb{R}^+$. The function $g(x, s)$ is measurable in $x$, of class $C^1$ in $s$ and it satisfies the following properties:

$$(14) \qquad \lim_{s \to +\infty} \frac{g(x, s)}{s^{q-1}} = 0, \qquad \lim_{s \to 0^+} \frac{g(x, s)}{s^{p-1}} = 0 \quad \text{uniformly in } x,$$

$$(15) \qquad 0 < \vartheta G(x,s) \leqslant g(x,s)s \quad \text{for } x \in \Lambda \text{ and } s \in \mathbb{R}^+,$$

$$(16) \qquad 0 \leqslant pG(x,s) \leqslant g(x,s)s \leqslant \frac{1}{\kappa} V(x)s^p \quad \text{for } x \in \Omega \setminus \Lambda \text{ and } s \in \mathbb{R}^+,$$

where we have set $G(x,s) := \int_0^s g(x,\tau)\,d\tau$.

Without loss of generality, we may assume that

$$(17) \qquad g(x,s) = 0 \quad \text{for a.e. } x \in \Omega \text{ and every } s < 0,$$

$$(18) \qquad j(x,s,\xi) = j(x,0,\xi) \quad \text{for every } x \in \Omega,\ s < 0 \text{ and } \xi \in \mathbb{R}^N.$$

Let now $J_\varepsilon : W_V(\Omega) \to \mathbb{R}$ be the functional defined as

$$J_\varepsilon(u) := \varepsilon^p \int_\Omega j(x,u,Du) + \frac{1}{p}\int_\Omega V(x)|u|^p - \int_\Omega G(x,u).$$

If $\bar{x}$ is in one of the $\Lambda_i$s, we also consider the "limit" functionals on $W^{1,p}(\mathbb{R}^N)$,

$$(19) \qquad I_{\bar{x}}(u) := \int_{\mathbb{R}^N} j(\bar{x},u,Du) + \frac{1}{p}\int_{\mathbb{R}^N} V(\bar{x})|u|^p - \int_{\mathbb{R}^N} F(u)$$

whose positive critical points solve equation (11). We denote by $c_{\bar{x}}$ the Mountain–Pass value of $I_{\bar{x}}$, namely

$$(20) \qquad c_{\bar{x}} := \inf_{\gamma \in \mathscr{P}_{\bar{x}}} \sup_{t \in [0,1]} I_{\bar{x}}(\gamma(t)),$$

$$(21) \qquad \mathscr{P}_{\bar{x}} := \Big\{ \gamma \in C([0,1], W^{1,p}(\mathbb{R}^N)) : \gamma(0) = 0,\ I_{\bar{x}}(\gamma(1)) < 0 \Big\}.$$

We set $c_i := c_{x_i}$ for every $i = 1, \ldots, k$. Considering $\sigma_i > 0$ such that

$$\sum_{i=1}^k \sigma_i < \frac{1}{2} \min\{c_i : i = 1, \ldots, k\},$$

we claim that, up to making $\Lambda_i$s smaller, we may assume that

$$(22) \qquad c_i \leqslant c_{\bar{x}} \leqslant c_i + \sigma_i \quad \text{for all } \bar{x} \in \Lambda_i.$$

In fact $c_i \leqslant c_{\bar{x}}$ follows because $x_i$ is a minimum of $V$ in $\Lambda_i$ and (12) holds. On the other hand, let us consider $\bar{x}_h \to x_i$ such that $\lim_h c_{\bar{x}_h} = \limsup_{\bar{x} \to x_i} c_{\bar{x}}$. Let $\gamma \in \mathscr{P}_{\bar{x}}$ be such that $\max_{\tau \in [0,1]} I_{x_i}(\gamma(\tau)) \leqslant c_i + \sigma_i$. Since $I_{\bar{x}_h} \to I_{x_i}$ uniformly on $\gamma$, we have that for $h$ large enough, $\gamma \in \mathscr{P}_{\bar{x}_h}$ and there exists $\tau_h \in [0,1]$ such that

$$c_{\bar{x}_h} \leqslant I_{\bar{x}_h}(\gamma(\tau_h)) \leqslant I_{x_i}(\gamma(\tau_h)) + o(1) \leqslant c_i + \sigma_i + o(1).$$

We deduce that $\limsup_{\bar{x} \to x_i} c_{\bar{x}} \leqslant c_i + \sigma_i$ so that the claim is proved.

If $\hat{\Lambda}_i$ denote mutually disjoint open sets compactly containing $\Lambda_i$, we introduce the functionals $J_{\varepsilon,i} : W^{1,p}(\hat{\Lambda}_i) \to \mathbb{R}$ as

$$(23) \qquad J_{\varepsilon,i}(u) := \varepsilon^p \int_{\hat{\Lambda}_i} j(x,u,Du) + \frac{1}{p}\int_{\hat{\Lambda}_i} V(x)|u|^p - \int_{\hat{\Lambda}_i} G(x,u)$$

for every $i = 1, \ldots, k$.

Finally, let us define the penalized functional $E_\varepsilon : W_V(\Omega) \to \mathbb{R}$ by setting

(24) $\qquad E_\varepsilon(u) := J_\varepsilon(u) + P_\varepsilon(u),$

$$P_\varepsilon(u) := M \sum_{i=1}^{k} \left( (J_{\varepsilon,i}(u)_+)^{1/2} - \varepsilon^{N/2}(c_i + \sigma_i)^{1/2} \right)_+^2,$$

where $M > 0$ is chosen so that

$$M > \frac{c_1 + \cdots + c_k}{\min\limits_{i=1,\ldots,k} \left\{ (2c_i)^{1/2} - (c_i + \sigma_i)^{1/2} \right\}}.$$

The functionals $J_\varepsilon$, $J_{\varepsilon,i}$ and $E_\varepsilon$ are merely continuous.

The next result provides the link between the critical points of $E_\varepsilon$ (see Definition 6.2) and the weak solutions of the original problem.

**Proposition 2.1.** *Let $u_\varepsilon \in W_V(\Omega)$ be any critical point of $E_\varepsilon$ and assume that there exists a positive number $\varepsilon_0$ such that the following conditions hold*

(25) $\qquad u_\varepsilon(x) < \ell \quad \text{for every } \varepsilon \in (0, \varepsilon_0) \text{ and } x \in \Omega \setminus \Lambda,$

(26) $\qquad \varepsilon^{-N} J_{\varepsilon,i}(u_\varepsilon) < c_i + \sigma_i \quad \text{for every } \varepsilon \in (0, \varepsilon_0) \text{ and } i = 1, \ldots, k.$

*Then, for every $\varepsilon \in (0, \varepsilon_0)$, $u_\varepsilon$ is a solution of $(P_\varepsilon)$.*

*Proof.* Let $\varepsilon \in (0, \varepsilon_0)$. By condition (26) and the definition of $P(u_\varepsilon)$, $u_\varepsilon$ is actually a critical point of $J_\varepsilon$. In view of (a) of Proposition 6.8, $u_\varepsilon$ is a weak solution of

$$-\varepsilon^p \text{div}\,(j_\xi(x, u, Du)) + \varepsilon^p j_s(x, u, Du) + V(x)|u|^{p-2}u = G(x, u).$$

Moreover, by (25) and the definition of $\widetilde{f}$, it results $G(x, u_\varepsilon(x)) = F(u_\varepsilon(x))$ for a.e. $x \in \Omega$. By (17) and (18) and arguing as in the proof of [28, Lemma 1], one gets $u_\varepsilon > 0$ in $\Omega$. Thus $u_\varepsilon$ is a solution of $(P_\varepsilon)$. □

The next Lemma is a variant of a local compactness property for bounded concrete Palais–Smale sequences (cf. [28, Theorem 2 and Lemma 3]; see also [5]).

**Lemma 2.2.** *Assume that (7), (8), (10) hold and let $(\psi_h) \subset L^\infty(\mathbb{R}^N)$ bounded with $\psi_h(x) \geqslant \lambda > 0$. Let $\varepsilon > 0$ and assume that $(u_h) \subset W^{1,p}(\mathbb{R}^N)$ is a bounded sequence such that*

$$\langle w_h, \varphi \rangle = \varepsilon^p \int_{\mathbb{R}^N} \psi_h(x) j_\xi(x, u_h, Du_h) \cdot D\varphi + \varepsilon^p \int_{\mathbb{R}^N} \psi_h(x) j_s(x, u_h, Du_h) \varphi$$

*for every $\varphi \in C_c^\infty(\mathbb{R}^N)$, where $(w_h)$ is strongly convergent in $W^{-1,p'}(\widetilde{\Omega})$ for a given bounded domain $\widetilde{\Omega}$ of $\mathbb{R}^N$.*

*Then $(u_h)$ admits a strongly convergent subsequence in $W^{1,p}(\widetilde{\Omega})$.*

Since $\Omega$ may be unbounded, in general the original functional $I_\varepsilon$ does not satisfy the concrete Palais–Smale condition (see Definition 6.6). In the following Lemma we prove that, instead, for every $\varepsilon > 0$ the functional $E_\varepsilon$ satisfies it at every level $c \in \mathbb{R}$.

**Lemma 2.3.** *Assume that conditions* (1), (5), (6), (7), (8), (9), (10) *hold. Let $\varepsilon > 0$.*

*Then $E_\varepsilon$ satisfies the concrete Palais–Smale condition at every level $c \in \mathbb{R}$.*

*Proof.* Let $(u_h) \subset W_V(\Omega)$ be a concrete Palais–Smale sequence for $E_\varepsilon$ at level $c$. We divide the proof into two steps:

<u>STEP I</u>. We prove that $(u_h)$ is bounded in $W_V(\Omega)$. From (15) and (16), we get

$$(27) \quad \vartheta \varepsilon^p \int_\Omega j(x, u_h, Du_h) + \frac{\vartheta}{p} \int_\Omega V(x)|u_h|^p \leqslant$$
$$\leqslant \int_\Lambda g(x, u_h) u_h + \frac{\vartheta}{p\kappa} \int_{\Omega \setminus \Lambda} V(x)|u_h|^p + \vartheta J_\varepsilon(u_h)$$

for every $h \in \mathbb{N}$. Moreover, by virtue of Proposition 6.4, for every $h \in \mathbb{N}$ we can compute $J'_\varepsilon(u_h)(u_h)$; in view of (16) we obtain

$$\int_\Lambda g(x, u_h) u_h + J'_\varepsilon(u_h)[u_h] \leqslant$$
$$\leqslant \varepsilon^p \int_\Omega j_\xi(x, u_h, Du_h) \cdot Du_h + \varepsilon^p \int_\Omega j_s(x, u_h, Du_h) u_h + \int_\Omega V(x)|u_h|^p$$

for every $h \in \mathbb{N}$. Notice that by (9) and the $p$–homogeneity of the map $\{\xi \mapsto j(x, s, \xi)\}$, it results

$$j_s(x, u_h, Du_h) u_h \leqslant \gamma j(x, u_h, Du_h),$$
$$j_\xi(x, u_h, Du_h) \cdot Du_h = p j(x, u_h, Du_h)$$

for every $h \in \mathbb{N}$. Therefore, we get

$$(28) \quad \int_\Lambda g(x, u_h) u_h + J'_\varepsilon(u_h)[u_h] \leqslant (\gamma + p)\varepsilon^p \int_\Omega j(x, u_h, Du_h) + \int_\Omega V(x)|u_h|^p$$

for every $h \in \mathbb{N}$. In view of (8), by combining inequalities (27) and (28) one gets

$$(29) \quad \min\left\{(\vartheta - \gamma - p)\nu\varepsilon^p, \frac{\vartheta}{p} - \frac{\vartheta}{p\kappa} - 1\right\} \int_\Omega \left(|Du_h|^p + V(x)|u_h|^p\right) \leqslant$$
$$\leqslant \vartheta J_\varepsilon(u_h) - J'_\varepsilon(u_h)[u_h]$$

for every $h \in \mathbb{N}$. In a similar fashion, arguing on the functionals $J_{\varepsilon,i}$, it results

$$(30) \quad \min\left\{(\vartheta - \gamma - p)\nu\varepsilon^p, \frac{\vartheta}{p} - \frac{\vartheta}{p\kappa} - 1\right\} \int_{\hat{\Lambda}_i} \left(|Du_h|^p + V(x)|u_h|^p\right) \leqslant$$
$$\leqslant \vartheta J_{\varepsilon,i}(u_h) - J'_{\varepsilon,i}(u_h)[u_h] \quad \text{for every } h \in \mathbb{N} \text{ and } i = 1, \ldots, k.$$

In particular, notice that one obtains

$$\bar{\vartheta} J_{\varepsilon,i}(u_h) - J'_{\varepsilon,i}(u_h)[u_h] \geqslant 0 \quad \text{for every } h \in \mathbb{N} \text{ and } i = 1, \ldots, k$$

and every $\gamma + p < \bar\vartheta < \vartheta$. Then, after some computations, one gets

$$\bar\vartheta P_\varepsilon(u_h) - P'_\varepsilon(u_h)[u_h] \geqslant$$

$$\geqslant -\bar\vartheta M \varepsilon^{N/2} \sum_{i=1}^{k} (c_i + \sigma_i)^{1/2} \left( (J_{\varepsilon,i}(u_h)_+)^{1/2} - \varepsilon^{N/2}(c_i + \sigma_i)^{1/2} \right)_+ \geqslant$$

$$\geqslant -C\varepsilon^{N/2} P_\varepsilon(u_h)^{1/2}$$

which implies, by Young's inequality, the existence of a constant $d > 0$ such that

(31) $$\vartheta P_\varepsilon(u_h) - P'_\varepsilon(u_h)[u_h] \geqslant -d\varepsilon^N$$

for every $h \in \mathbb{N}$. By combining (29) with (31), since

$$E_\varepsilon(u_h) = c + o(1), \quad E'_\varepsilon(u_h)[u_h] = o(\|u_h\|_{W_V})$$

as $h \to +\infty$, one obtains

(32) $$\int_\Omega \left( |Du_h|^p + V(x)|u_h|^p \right) \leqslant$$

$$\leqslant \frac{\vartheta c + d\varepsilon^N}{\min\left\{ (\vartheta - \gamma - p)\nu\varepsilon^p, \frac{\vartheta}{p} - \frac{\vartheta}{p\kappa} - 1 \right\}} + o(\|u_h\|_{W_V}) + o(1)$$

as $h \to +\infty$, which yields the boundedness of $(u_h)$ in $W_V(\Omega)$.

<u>Step II</u>. By virtue of Step I, there exists $u \in W_V(\Omega)$ such that, up to a subsequence, $(u_h)$ weakly converges to $u$ in $W_V(\Omega)$. Let us now prove that actually $(u_h)$ converges strongly to $u$ in $W_V(\Omega)$. If we define for every $h \in \mathbb{N}$ the weights

$$\theta_{h,i} = M\left[ (J_{\varepsilon,i}(u_h)_+)^{1/2} - \varepsilon^{N/2}(c_i + \sigma_i)^{1/2} \right]_+ (J_{\varepsilon,i}(u_h)_+)^{-1/2}, \quad i = 1, \ldots, k$$

and put $\theta_h(x) = \sum_{i=1}^{k} \theta_{h,i} \mathbb{1}_{\hat\Lambda_i}(x)$ with $0 \leqslant \theta_{h,i} \leqslant M$, after a few computations, one gets

$$\langle w_h, \varphi \rangle = \varepsilon^p \int_\Omega (1 + \theta_h) j_\xi(x, u_h, Du_h) \cdot D\varphi + \varepsilon^p \int_\Omega (1 + \theta_h) j_s(x, u_h, Du_h)\varphi$$

for every $\varphi \in C_c^\infty(\Omega)$, where

$$w_h = (1 + \theta_h)g(x, u_h) - (1 + \theta_h)V(x)|u_h|^{p-2}u_h + \xi_h,$$

with $\xi_h \to 0$ strongly in $W^{-1,p'}(\Omega)$. Since, up to a subsequence, $(w_h)$ strongly converges to $w := (1+\bar\theta)g(x,u) - (1+\bar\theta)V(x)|u|^{p-2}u$ in $W^{-1,p'}(B_\varrho)$ for every $\varrho > 0$, by applying Lemma 2.2 with $\widetilde\Omega = B_\varrho \cap \Omega$ and $\psi_h(x) = 1 + \theta_h(x)$, it suffices to show that for every $\delta > 0$ there exists $\varrho > 0$ such that

(33) $$\limsup_h \int_{\Omega \setminus B_\varrho} \left( |Du_h|^p + V(x)|u_h|^p \right) < \delta.$$

Consider a cut–off function $\chi_\varrho \in C^\infty(\mathbb{R}^N)$ with $0 \leqslant \chi_\varrho \leqslant 1$, $\chi_\varrho = 0$ on $B_{\varrho/2}$, $\chi_\varrho = 1$ on $\mathbb{R}^N \setminus B_\varrho$ and $|D\chi_\varrho| \leqslant a/\varrho$ for some $a > 0$. By taking $\varrho$ large enough, we have

$$(34) \qquad \bigcup_{i=1}^{k} \hat{\Lambda}_i \cap \mathrm{supt}(\chi_\varrho) = \emptyset.$$

Let now $\zeta : \mathbb{R} \to \mathbb{R}$ be the map defined by

$$(35) \qquad \zeta(s) := \begin{cases} 0 & \text{if } s < 0 \\ \bar{M}s & \text{if } 0 \leqslant s < R \\ \bar{M}R & \text{if } s \geqslant R, \end{cases}$$

being $R > 0$ the constant defined in (10) and $\bar{M}$ a positive number (which exists by the growths (7) and (8)) such that

$$(36) \qquad |j_s(x,s,\xi)| \leqslant p\bar{M} j(x,s,\xi)$$

for every $x \in \Omega$, $s \in \mathbb{R}$ and $\xi \in \mathbb{R}^N$. Notice that, by combining (10) and (36), we obtain

$$(37) \quad j_s(x,s,\xi) + p\zeta'(s)j(x,s,\xi) \geqslant 0 \quad \text{for every } x \in \Omega,\ s \in \mathbb{R} \text{ and } \xi \in \mathbb{R}^N.$$

By (34) it is easily proved that $P'_\varepsilon(u_h)(\chi_\varrho u_h e^{\zeta(u_h)}) = 0$ for every $h$. Therefore, since the sequence $(\chi_\varrho u_h e^{\zeta(u_h)})$ is bounded in $W_V(\Omega)$, taking into account (37) and (18) we obtain

$$\begin{aligned}
o(1) = J'_\varepsilon(u_h)(\chi_\varrho u_h e^{\zeta(u_h)}) &= \\
&= \varepsilon^p \int_\Omega j_\xi(x, u_h, Du_h) \cdot Du_h \chi_\varrho e^{\zeta(u_h)} + \\
&\quad + \varepsilon^p \int_\Omega j_\xi(x, u_h, Du_h) \cdot D\chi_\varrho u_h e^{\zeta(u_h)} + \\
&\quad + \varepsilon^p \int_\Omega \left[ j_s(x, u_h, Du_h) + p\zeta'(u_h)j(x, u_h, Du_h) \right] u_h \chi_\varrho e^{\zeta(u_h)} + \\
&\quad + \int_\Omega V(x)|u_h|^p \chi_\varrho e^{\zeta(u_h)} - \int_\Omega g(x, u_h) u_h \chi_\varrho e^{\zeta(u_h)} \geqslant \\
&\geqslant \int_\Omega \left( p\varepsilon^p j(x, u_h, Du_h) + V(x)|u_h|^p \right) \chi_\varrho e^{\zeta(u_h)} + \\
&\quad + \varepsilon^p \int_\Omega j_\xi(x, u_h, Du_h) \cdot D\chi_\varrho u_h e^{\zeta(u_h)} + \\
&\quad - \int_\Omega g(x, u_h) u_h \chi_\varrho e^{\zeta(u_h)}
\end{aligned}$$

as $h \to +\infty$. Therefore, in view of (16) and (34), it results

$$o(1) \geqslant \int_\Omega \left(p\varepsilon^p \nu |Du_h|^p + V(x)|u_h|^p\right)\chi_\varrho e^{\zeta(u_h)} +$$

$$+ \varepsilon^p \int_\Omega j_\xi(x, u_h, Du_h) \cdot D\chi_\varrho u_h e^{\zeta(u_h)} - \frac{1}{\kappa}\int_\Omega V(x)|u_h|^p \chi_\varrho e^{\zeta(u_h)}$$

as $h \to +\infty$ for $\varrho$ large enough. Since by (7) we have

$$\left|\int_\Omega j_\xi(x, u_h, Du_h) \cdot D\chi_\varrho u_h e^{\zeta(u_h)}\right| \leqslant \frac{C}{\varrho}\|Du_h\|_p^{p-1}\|u_h\|_p \leqslant \frac{\widetilde{C}}{\varrho},$$

there exists a positive constant $C'$ such that

$$\limsup_h \int_{\Omega \setminus B_\varrho} \left(|Du_h|^p + V(x)|u_h|^p\right) \leqslant \frac{C'}{\varrho}$$

which yields (33). The proof is now complete. □

3. Two consequences of a Pucci–Serrin type identity

Let $\mathscr{L} : \mathbb{R}^N \times \mathbb{R} \times \mathbb{R}^N \to \mathbb{R}$ be a function of class $C^1$ such that the function $\{\xi \mapsto \mathscr{L}(x, s, \xi)\}$ is strictly convex for every $(x, s) \in \mathbb{R}^N \times \mathbb{R}$, and let $\varphi \in L^\infty_{\text{loc}}(\mathbb{R}^N)$.

We now recall a Pucci–Serrin variational identity for locally Lipschitz continuous solutions of a general class of Euler equations, recently obtained in [10]. Notice that the classical identity [24] is not applicable here, since it requires the $C^2$ regularity of the solutions while in our degenerate setting the maximal regularity is $C^{1,\beta}_{\text{loc}}$ (see [14, 31]).

**Theorem 3.1.** *Let $u : \mathbb{R}^N \to \mathbb{R}$ be a locally Lipschitz solution of*

$$-\text{div}\,(D_\xi \mathscr{L}(x, u, Du)) + D_s \mathscr{L}(x, u, Du) = \varphi \quad \text{in } \mathscr{D}'(\mathbb{R}^N).$$

*Then*

(38)
$$\sum_{i,j=1}^N \int_{\mathbb{R}^N} D_i h^j D_{\xi_i}\mathscr{L}(x, u, Du) D_j u +$$
$$- \int_{\mathbb{R}^N} \left[(\text{div}\,h)\,\mathscr{L}(x, u, Du) + h \cdot D_x \mathscr{L}(x, u, Du)\right] = \int_{\mathbb{R}^N}(h \cdot Du)\varphi$$

*for every $h \in C^1_c(\mathbb{R}^N, \mathbb{R}^N)$.*

We want to derive two important consequences of the previous variational identity.

In the first we show that the Mountain–Pass value associated with a large class of elliptic autonomous equations is the minimal among other nontrivial critical values.

**Theorem 3.2.** *Let $\bar{x} \in \mathbb{R}^N$ and assume that conditions (1), (5), (6), (7), (8), (9), (10) hold. Then the equation*

(39) $\quad -\text{div}\,(j_\xi(\bar{x}, u, Du)) + j_s(\bar{x}, u, Du) + V(\bar{x})u^{p-1} = f(u) \quad \text{in } \mathbb{R}^N$

admits a least energy solution $u \in W^{1,p}(\mathbb{R}^N)$, that is

$$I_{\bar{x}}(u) = \inf \{I_{\bar{x}}(w) : w \in W^{1,p}(\mathbb{R}^N) \setminus \{0\} \text{ is a solution of } (39)\},$$

where $I_{\bar{x}}$ is as in (20). Moreover, $I_{\bar{x}}(u) = c_{\bar{x}}$, that is $u$ is at the Mountain–Pass level.

*Proof.* We divide the proof into two steps.

STEP I. Let $u$ be any nontrivial solution of (39), and let us prove that $I_{\bar{x}}(u) \geq c_{\bar{x}}$. By the assumptions on $V$ and $f$, it is readily seen that there exist $\varrho_0 > 0$ and $\delta_0 > 0$ such that $I_{\bar{x}}(v) \geq \delta_0$ for every $v \in W^{1,p}(\mathbb{R}^N)$ with $\|v\|_{1,p} = \varrho_0$. In particular $I_{\bar{x}}$ has a Mountain–Pass geometry. As we will see, $\mathscr{P}_{\bar{x}} \neq \emptyset$, so that $c_{\bar{x}}$ is well defined. Let now $u$ be a positive solution of (39) and consider the dilation path

$$\gamma(t)(x) := \begin{cases} u\left(\frac{x}{t}\right) & \text{if } t > 0 \\ 0 & \text{if } t = 0. \end{cases}$$

Notice that $\|\gamma(t)\|_{1,p}^p = t^{N-p}\|Du\|_p^p + t^N \|u\|_p^p$ for every $t \in \mathbb{R}^+$, which implies that the curve $\gamma$ belongs to $C(\mathbb{R}^+, W^{1,p}(\mathbb{R}^N))$. For the sake of simplicity, we consider the continuous function $H : \mathbb{R}^+ \to \mathbb{R}$ defined by

$$H(s) = \int_0^s h(t)\, dt, \quad \text{where } h(s) = -V(\bar{x})s^{p-1} + f(s).$$

For every $t \in \mathbb{R}^+$ it results that

$$I_{\bar{x}}(\gamma(t)) = \int_{\mathbb{R}^N} j(\bar{x}, \gamma(t), D\gamma(t)) - \int_{\mathbb{R}^N} H(\gamma(t)) =$$
$$= t^{N-p} \int_{\mathbb{R}^N} j(\bar{x}, u, Du) - t^N \int_{\mathbb{R}^N} H(u)$$

which yields, for every $t \in \mathbb{R}^+$

$$(40) \quad \frac{d}{dt} I_{\bar{x}}(\gamma(t)) = (N-p)t^{N-p-1} \int_{\mathbb{R}^N} j(\bar{x}, u, Du) - Nt^{N-1} \int_{\mathbb{R}^N} H(u).$$

By virtue of (8) and (10), a standard argument yields $u \in L^\infty_{\text{loc}}(\mathbb{R}^N)$ (see [26, Theorem 1]); by the regularity results of [14, 31], it follows that $u \in C^{1,\beta}_{\text{loc}}(\mathbb{R}^N)$ for some $0 < \beta < 1$. Then, since $\{\xi \mapsto j(x, s, \xi)\}$ is strictly convex, we can use Theorem 3.1 by choosing in (38) $\varphi = 0$ and

$$(41) \quad \mathscr{L}(s, \xi) := j(\bar{x}, s, \xi) - H(s) \quad \text{for every } s \in \mathbb{R}^+ \text{ and } \xi \in \mathbb{R}^N,$$
$$h(x) = h_k(x) := T\left(\frac{x}{k}\right) x \quad \text{for every } x \in \mathbb{R}^N \text{ and } k \geq 1,$$

being $T \in C^1_c(\mathbb{R}^N)$ such that $T(x) = 1$ if $|x| \leq 1$ and $T(x) = 0$ if $|x| \geq 2$. In particular, for every $k$ we have that $h_k \in C^1_c(\mathbb{R}^N, \mathbb{R}^N)$ and

$$D_i h_k^j(x) = D_i T\left(\frac{x}{k}\right) \frac{x_j}{k} + T\left(\frac{x}{k}\right) \delta_{ij} \quad \text{for every } x \in \mathbb{R}^N, i,j = 1,\ldots,N,$$
$$(\text{div } h_k)(x) = DT\left(\frac{x}{k}\right) \cdot \frac{x}{k} + NT\left(\frac{x}{k}\right) \quad \text{for every } x \in \mathbb{R}^N.$$

Then, since $D_x \mathscr{L}(u, Du) = 0$, it follows by (38) that

$$\sum_{i,j=1}^{n} \int_{\mathbb{R}^N} D_i T\left(\frac{x}{k}\right) \frac{x_j}{k} D_j u D_{\xi_i} \mathscr{L}(u, Du) + \int_{\mathbb{R}^N} T\left(\frac{x}{k}\right) D_\xi \mathscr{L}(u, Du) \cdot Du +$$

$$- \int_{\mathbb{R}^N} DT\left(\frac{x}{k}\right) \cdot \frac{x}{k} \mathscr{L}(u, Du) - \int_{\mathbb{R}^N} NT\left(\frac{x}{k}\right) \mathscr{L}(u, Du) = 0$$

for every $k \geqslant 1$. Since there exists $C > 0$ with

$$D_i T\left(\frac{x}{k}\right) \frac{x_j}{k} \leqslant C \quad \text{for every } x \in \mathbb{R}^N,\ k \geqslant 1 \text{ and } i, j = 1, \ldots, N,$$

by the Dominated Convergence Theorem, letting $k \to +\infty$, we obtain

$$\int_{\mathbb{R}^N} \left[ N\mathscr{L}(u, Du) - D_\xi \mathscr{L}(u, Du) \cdot Du \right] = 0,$$

namely, by (41) and the $p$–homogeneity of $\{\xi \mapsto j(x, s, \xi)\}$,

(42) $$(N - p) \int_{\mathbb{R}^N} j(\bar{x}, u, Du) = N \int_{\mathbb{R}^N} H(u).$$

In particular notice that $\int_{\mathbb{R}^N} H(u) > 0$. By plugging this formula into (40), we obtain

$$\frac{d}{dt} I_{\bar{x}}(\gamma(t)) = N(1 - t^p) t^{N-p-1} \int_{\mathbb{R}^N} H(u)$$

which yields $\frac{d}{dt} I_{\bar{x}}(\gamma(t)) > 0$ for $0 < t < 1$ and $\frac{d}{dt} I_{\bar{x}}(\gamma(t)) < 0$ for $t > 1$, namely

$$\sup_{t \in [0,+\infty[} I_{\bar{x}}(\gamma(t)) = I_{\bar{x}}(\gamma(1)) = I_{\bar{x}}(u).$$

Moreover, observe that $\gamma(0) = 0$ and $I_{\bar{x}}(\gamma(T)) < 0$ for $T > 0$ sufficiently large. Then, after a suitable scale change in $t$, $\gamma \in \mathscr{P}_{\bar{x}}$ and the assertion follows.

<u>STEP II</u> Let us now prove that (39) has a nontrivial solution $u \in W^{1,p}(\mathbb{R}^N)$ such that $c_{\bar{x}} \geqslant I_{\bar{x}}(u)$. Let $(u_h)$ be a Palais–Smale sequence for $I_{\bar{x}}$ at the level $c_{\bar{x}}$. Since $(u_h)$ is bounded in $W^{1,p}(\mathbb{R}^N)$, considering the test $u_h e^{\zeta(u_h)}$ with $\zeta$ as in (35), and recalling (37), we have

$$pc_{\bar{x}} + o(1) = pI_{\bar{x}}(u_h) - I'_{\bar{x}}(u_h)[u_h e^{\zeta(u_h)}] =$$

$$= \int_{\mathbb{R}^N} p(1 - e^{\zeta(u_h)}) j(\bar{x}, u_h, Du_h) + \int_{\mathbb{R}^N} (1 - e^{\zeta(u_h)}) V(\bar{x}) |u_h|^p +$$

$$- \int_{\mathbb{R}^N} \left[ p\zeta'(u_h) j(\bar{x}, u_h, Du_h) + j_s(\bar{x}, u_h, Du_h) \right] u_h e^{\zeta(u_h)} +$$

$$- \int_{\mathbb{R}^N} pF(u_h) + \int_{\mathbb{R}^N} f(u_h) u_h e^{\zeta(u_h)} \leqslant$$

$$\leqslant - \int_{\mathbb{R}^N} pF(u_h) + \int_{\mathbb{R}^N} f(u_h) u_h e^{\zeta(u_h)} \leqslant C \int_{\mathbb{R}^N} |u_h|^p + |u_h|^q$$

for some $C > 0$. By [21, Lemma I.1], we conclude that $(u_h)$ may not vanish in $L^p$, that is there exists $x_h \in \mathbb{R}^N$, $R > 0$ and $\lambda > 0$ such that for $h$ large

$$(43) \qquad \int_{x_h + B_R} |u_h|^p \geqslant \lambda.$$

Let $v_h(x) := u_h(x_h + x)$ and let $u \in W^{1,p}(\mathbb{R}^N)$ be such that $v_h \rightharpoonup u$ weakly in $W^{1,p}(\mathbb{R}^N)$. Since $v_h$ is a Palais–Smale sequence for $I_{\bar{x}}$ at level $c_{\bar{x}}$, by Lemma 2.2, we have that $v_h \to u$ strongly in $W^{1,p}_{\mathrm{loc}}(\mathbb{R}^N)$. By (43), we deduce that $u$ is a nontrivial solution of (39). Let $\delta > 0$; we claim that there exists $\varrho > 0$ such that

$$(44) \qquad \liminf_h \int_{\mathbb{R}^N \setminus B_\varrho} \left[ j(\bar{x}, v_h, Dv_h) + \frac{1}{p} V(\bar{x}) |v_h|^p - F(v_h) \right] \geqslant -\delta.$$

In fact, let $\varrho > 0$, and let $\eta_\varrho$ be a smooth function such that $0 \leqslant \eta_\varrho \leqslant 1$, $\eta_\varrho = 0$ on $B_{\varrho-1}$, $\eta_\varrho = 1$ on $\mathbb{R}^N \setminus B_\varrho$ and $\|D\eta_\varrho\|_\infty \leqslant 2$. By Proposition 6.4, testing with $\eta_\varrho v_h$, we get

$$\langle w_h, \eta_\varrho v_h \rangle - \int_{B_\varrho \setminus B_{\varrho-1}} \left[ j_\xi(\bar{x}, v_h, Dv_h) \cdot D(\eta_\varrho v_h) + j_s(\bar{x}, v_h, Dv_h) \eta_\varrho v_h + \right.$$
$$+ V(\bar{x}) |v_h|^p \eta_\varrho - f(v_h) v_h \eta_\varrho \Big] =$$
$$= \int_{\mathbb{R}^N \setminus B_\varrho} \left[ j_\xi(\bar{x}, v_h, Dv_h) \cdot D(\eta_\varrho v_h) + j_s(\bar{x}, v_h, Dv_h) \eta_\varrho v_h + \right.$$
$$+ V(\bar{x}) |v_h|^p \eta_\varrho - f(v_h) v_h \eta_\varrho \Big]$$

where $w_h \to 0$ strongly in $W^{-1,p'}(\mathbb{R}^N)$. For the right hand side we have

$$\int_{\mathbb{R}^N \setminus B_\varrho} \left[ j_\xi(\bar{x}, v_h, Dv_h) \cdot D(\eta_\varrho v_h) + j_s(\bar{x}, v_h, Dv_h) \eta_\varrho v_h + \right.$$
$$+ V(\bar{x}) |v_h|^p \eta_\varrho - f(v_h) v_h \eta_\varrho \Big] =$$
$$= \int_{\mathbb{R}^N \setminus B_\varrho} \left[ pj(\bar{x}, v_h, Dv_h) + j_s(\bar{x}, v_h, Dv_h) v_h + \right.$$
$$+ V(\bar{x}) |v_h|^p - f(v_h) v_h \Big],$$

and by (9) we have

$$\int_{\mathbb{R}^N \setminus B_\varrho} \left[ pj(\bar{x}, v_h, Dv_h) + j_s(\bar{x}, v_h, Dv_h)v_h + V(\bar{x})|v_h|^p - f(v_h)v_h \right] \leqslant$$

$$\leqslant (p+\gamma) \int_{\mathbb{R}^N \setminus B_\varrho} j(\bar{x}, v_h, Dv_h) + \int_{\mathbb{R}^N \setminus B_\varrho} V(\bar{x})|v_h|^p - f(v_h)v_h =$$

$$= (p+\gamma) \int_{\mathbb{R}^N \setminus B_\varrho} \left[ j(\bar{x}, v_h, Dv_h) + \frac{1}{p}V(\bar{x})|v_h|^p - F(v_h) \right] +$$

$$- \frac{p+\gamma}{p} \int_{\mathbb{R}^N \setminus B_\varrho} V(\bar{x})|v_h|^p + \int_{\mathbb{R}^N \setminus B_\varrho} V(\bar{x})|v_h|^p +$$

$$+ \int_{\mathbb{R}^N \setminus B_\varrho} \left[ (p+\gamma)F(v_h) - f(v_h)v_h \right] \leqslant$$

$$\leqslant (p+\gamma) \int_{\mathbb{R}^N \setminus B_\varrho} \left[ j(\bar{x}, v_h, Dv_h) + \frac{1}{p}V(\bar{x})|v_h|^p - F(v_h) \right] +$$

$$+ \int_{\mathbb{R}^N \setminus B_\varrho} \left[ (p+\gamma)F(v_h) - \vartheta F(v_h) \right] \leqslant$$

$$\leqslant (p+\gamma) \int_{\mathbb{R}^N \setminus B_\varrho} \left[ j(\bar{x}, v_h, Dv_h) + \frac{1}{p}V(\bar{x})|v_h|^p - F(v_h) \right].$$

We conclude that

$$(p+\gamma) \int_{\mathbb{R}^N \setminus B_\varrho} \left[ j(\bar{x}, v_h, Dv_h) + \frac{1}{p}V(\bar{x})|v_h|^p - F(v_h) \right] \geqslant \langle w_h, \eta_\varrho v_h \rangle +$$

$$- \int_{B_\varrho \setminus B_{\varrho-1}} \left[ j_\xi(\bar{x}, v_h, Dv_h) \cdot D(\eta_\varrho v_h) + j_s(\bar{x}, v_h, Dv_h)\eta_\varrho v_h + \right.$$

$$\left. + V(\bar{x})|v_h|^p \eta_\varrho - f(v_h)v_h \eta_\varrho \right].$$

Since by Lemma 2.2 we have $v_h \to u$ strongly in $W^{1,p}(B_\varrho)$, we get

$$\lim_h \int_{B_\varrho \setminus B_{\varrho-1}} \left[ j_\xi(\bar{x}, v_h, Dv_h) \cdot D(\eta_\varrho v_h) + j_s(\bar{x}, v_h, Dv_h)\eta_\varrho v_h + \right.$$

$$\left. + V(\bar{x})|v_h|^p \eta_\varrho - f(v_h)v_h \eta_\varrho \right] =$$

$$= \int_{B_\varrho \setminus B_{\varrho-1}} \left[ j_\xi(\bar{x}, u, Du) \cdot D(\eta_\varrho u) + j_s(\bar{x}, u, Du)\eta_\varrho u + \right.$$

$$\left. + V(\bar{x})|u|^p \eta_\varrho - f(u)u\eta_\varrho \right],$$

and so we deduce that for every $\delta > 0$ there exists $\bar{\varrho} > 0$ such that for all $\varrho > \bar{\varrho}$ we have

$$\liminf_h \int_{\mathbb{R}^N \setminus B_\varrho} \left[ j(\bar{x}, v_h, Dv_h) + \frac{1}{p}V(\bar{x})|v_h|^p - F(v_h) \right] \geqslant -\delta.$$

Futhermore we have
$$\lim_h \int_{B_\varrho} \left[j(\bar{x}, v_h, Dv_h) + \frac{1}{p}V(\bar{x})|v_h|^p - F(v_h)\right] = I_{\bar{x}}(u, B_\varrho),$$
where
$$I_{\bar{x}}(u, B_\varrho) := \int_{B_\varrho} \left[j(\bar{x}, u, Du) + \frac{1}{p}V(\bar{x})|u|^p - F(u)\right],$$
and so we conclude that for all $\varrho > \bar{\varrho}$
$$c_{\bar{x}} \geqslant I_{\bar{x}}(u, B_\varrho) - \delta.$$
Letting $\varrho \to +\infty$ and since $\delta$ is arbitrary, we get
$$c_{\bar{x}} \geqslant I_{\bar{x}}(u),$$
and the proof is concluded. □

The second result can be considered as an extension (also with a different proof) of [12, Lemma 2.3] to a general class of elliptic equations. Again we stress that, in this degenerate setting, Theorem 3.1 plays an important role.

**Lemma 3.3.** *Let $u \in W^{1,p}(\mathbb{R}^N)$ be a positive solution of the equation*
$$(45) \quad -\mathrm{div}\,(j_\xi(\bar{x}, u, Du)) + j_s(\bar{x}, u, Du) + V(\bar{x})u^{p-1} =$$
$$= \mathbb{1}_{\{x_1 < 0\}}(x)f(u) + \mathbb{1}_{\{x_1 > 0\}}(x)\widetilde{f}(u) \quad \text{in } \mathbb{R}^N.$$

*Then $u$ is actually a solution of the equation*
$$(46) \quad -\mathrm{div}\,(j_\xi(\bar{x}, u, Du)) + j_s(\bar{x}, u, Du) + V(\bar{x})u^{p-1} = f(u) \quad \text{in } \mathbb{R}^N.$$

*Proof.* Let us first show that $u(x) \leqslant \ell$ on the set $\{x_1 = 0\}$. As in the proof of Theorem 3.2 it follows that $u \in C^{1,\beta}_{\mathrm{loc}}(\mathbb{R}^N)$ for some $0 < \beta < 1$. Then we can apply again Theorem 3.1 by choosing this time in (38):
$$\mathscr{L}(s, \xi) := j(\bar{x}, s, \xi) + \frac{V(\bar{x})}{p}s^p \quad \text{for every } s \in \mathbb{R}^+ \text{ and } \xi \in \mathbb{R}^N,$$
$$\varphi(x) := \mathbb{1}_{\{x_1 < 0\}}(x)f(u(x)) + \mathbb{1}_{\{x_1 > 0\}}(x)\widetilde{f}(u(x)) \quad \text{for every } x \in \mathbb{R}^N,$$
$$h(x) = h_k(x) := \left(T\left(\frac{x}{k}\right), 0, \ldots, 0\right) \quad \text{for every } x \in \mathbb{R}^N \text{ and } k \geqslant 1$$
being $T \in C^1_c(\mathbb{R}^N)$ such that $T(x) = 1$ if $|x| \leqslant 1$ and $T(x) = 0$ if $|x| \geqslant 2$. Then $h_k \in C^1_c(\mathbb{R}^N, \mathbb{R}^N)$ and, taking into account that $D_x\mathscr{L}(u, Du) = 0$, we have
$$\int_{\mathbb{R}^N} \left[\frac{1}{k}\sum_{i=1}^N D_iT\left(\frac{x}{k}\right)D_1u D_{\xi_i}\mathscr{L}(u, Du) - D_1T\left(\frac{x}{k}\right)\frac{1}{k}\mathscr{L}(u, Du)\right] =$$
$$= \int_{\mathbb{R}^N} T\left(\frac{x}{k}\right)\varphi(x)D_1u$$

for every $k \geqslant 1$. Again by the Dominated Convergence Theorem, letting $k \to +\infty$, it follows $\int_{\mathbb{R}^N} \varphi(x) D_{x_1} u = 0$, that is, after integration by parts,

$$\int_{\mathbb{R}^{N-1}} \left[ F(u(0, x')) - \widetilde{F}(u(0, x')) \right] dx' = 0.$$

Taking into account that $F(s) \geqslant \widetilde{F}(s)$ with equality only if $s \leqslant \ell$, we get

(47) $$u(0, x') \leqslant \ell \quad \text{for every } x' \in \mathbb{R}^{N-1}.$$

To prove that actually

(48) $$u(x_1, x') \leqslant \ell \quad \text{for every } x_1 > 0 \text{ and } x' \in \mathbb{R}^{N-1},$$

let us test equation (45) with the function

$$\eta(x) = \begin{cases} 0 & \text{if } x_1 < 0 \\ (u(x_1, x') - \ell)^+ e^{\zeta(u(x_1, x'))} & \text{if } x_1 > 0, \end{cases}$$

where $\zeta : \mathbb{R}^+ \to \mathbb{R}$ is the map defined in (35). Notice that, in view of (47), the function $\eta$ belongs to $W^{1,p}(\mathbb{R}^N)$. After some computations, one obtains

(49) $$\int_{\{x_1 > 0\}} pj(\bar{x}, u, D(u-\ell)^+) e^{\zeta(u)} +$$
$$+ \int_{\{x_1 > 0\}} \left[ j_s(\bar{x}, u, Du) + p\zeta'(u) j(\bar{x}, u, Du) \right] (u-\ell)^+ e^{\zeta(u)} +$$
$$+ \int_{\{x_1 > 0\}} \left[ V(\bar{x}) - \frac{\alpha}{\kappa} \right] u^{p-1} (u-\ell)^+ e^{\zeta(u)} = 0.$$

By (1) and (37) all the terms in (49) must be equal to zero. We conclude that $(u - \ell)^+ = 0$ on $\{x_1 > 0\}$, namely (48) holds. In particular $\varphi(x) = f(u(x))$ for every $x \in \mathbb{R}^N$, so that $u$ is a solution of (46). □

## 4. Energy estimates

Let $d_{\varepsilon,i}$ be the Mountain–Pass critical value which corresponds to the functional $J_{\varepsilon,i}$ defined in (23). More precisely,

(50) $$d_{\varepsilon,i} := \inf_{\gamma_i \in \Gamma_i} \sup_{t \in [0,1]} J_{\varepsilon,i}(\gamma_i(t))$$

where

$$\Gamma_i := \left\{ \gamma_i \in C([0,1], W^{1,p}(\widehat{\Lambda}_i)) : \gamma_i(0) = 0, J_{\varepsilon,i}(\gamma_i(1)) < 0 \right\}.$$

Then the following result holds.

**Lemma 4.1.** *We have*

$$\lim_{\varepsilon \to 0^+} \varepsilon^{-N} d_{\varepsilon,i} = c_i$$

*for every $i = 1, \ldots, k$.*

*Proof.* The inequality

$$d_{\varepsilon,i} \leqslant \varepsilon^N c_i + o(\varepsilon^N) \tag{51}$$

can be easily derived (see the first part of the proof of Lemma 4.2). Let us prove the opposite inequality, which is harder. To this aim, we divide the proof into two steps.

<u>Step I</u>. Let $w_\varepsilon$ be a Mountain–Pass critical point for $J_{\varepsilon,i}$. We have $w_\varepsilon \geqslant 0$, and by regularity results $w_\varepsilon \in L^\infty(\hat{\Lambda}_i) \cap C^{1,\alpha}_{\text{loc}}(\hat{\Lambda}_i)$. Let us define

$$M_\varepsilon := \sup_{x \in \hat{\Lambda}_i} w_\varepsilon(x) < +\infty,$$

and for all $\delta > 0$ define the set

$$U_\delta := \{x \in \hat{\Lambda}_i : w_\varepsilon(x) > M_\varepsilon - \delta\}.$$

By Proposition 6.4, we may use the following nontrivial test for the equation satisfied by $w_\varepsilon$

$$\varphi_\delta := [w_\varepsilon - (M_\varepsilon - \delta)]^+ e^{\zeta(w_\varepsilon)},$$

where the map $\zeta : \mathbb{R}^+ \to \mathbb{R}$ is defined as in (35). We have

$$D\varphi_\delta = e^{\zeta(w_\varepsilon)} Dw_\varepsilon \mathbb{1}_{U_\delta} + \varphi_\delta \zeta'(w_\varepsilon) Dw_\varepsilon,$$

and so we obtain

$$\varepsilon^p \int_{U_\delta} p j(x, w_\varepsilon, Dw_\varepsilon) e^{\zeta(w_\varepsilon)} +$$
$$+ \varepsilon^p \int_{U_\delta} \left[p\zeta'(w_\varepsilon) j(x, w_\varepsilon, Dw_\varepsilon) + j_s(x, w_\varepsilon, Dw_\varepsilon)\right] \varphi_\delta =$$
$$= \int_{U_\delta} \left[-V(x) w_\varepsilon^{p-1} + g(x, w_\varepsilon)\right] \varphi_\delta.$$

Then, by (37), it results

$$\int_{U_\delta} \left[-V(x) w_\varepsilon^{p-1} + g(x, w_\varepsilon)\right] \varphi_\delta \geqslant \varepsilon^p \int_{U_\delta} p j(x, w_\varepsilon, Dw_\varepsilon) e^{\zeta(w_\varepsilon)} > 0. \tag{52}$$

Suppose that $U_\delta \cap \Lambda_i = \emptyset$ for some $\delta > 0$; we have that $g(x, w_\varepsilon) = \widetilde{f}(w_\varepsilon)$ on $U_\delta$, so that

$$\int_{U_\delta} \left[-V(x) w_\varepsilon^{p-1} + \widetilde{f}(w_\varepsilon)\right] \varphi_\delta > 0. \tag{53}$$

On the other hand, we note that by construction $\widetilde{f}(w_\varepsilon) \leqslant \frac{1}{k} V(x) w_\varepsilon^{p-1}$ with strict inequality on an open subset of $U_\delta$. We deduce that (53) cannot hold, and so $U_\delta \cap \Lambda_i \neq \emptyset$ for all $\delta$. Since $\Lambda_i$ is compact, we conclude that $w_\varepsilon$ admits a maximum point $x_\varepsilon$ in $\Lambda_i$. Moreover, we have $w_\varepsilon(x_\varepsilon) \geqslant \ell$, where $\ell$ is as in (13), since otherwise (52) cannot hold.

Let us now consider the functions $v_\varepsilon(y) := w_\varepsilon(x_\varepsilon + \varepsilon y)$ and let $\varepsilon_j \to 0$. We have that, up to a subsequence, $x_{\varepsilon_j} \to \bar{x} \in \Lambda_i$. Since $w_\varepsilon$ is a Mountain–Pass critical point of $J_{\varepsilon,i}$, arguing as in Step I of Lemma 2.3 there exists $C > 0$ such that
$$\int_{\mathbb{R}^N} \left( \varepsilon^p |Dw_\varepsilon|^p + V(x)|w_\varepsilon|^p \right) \leqslant C d_{\varepsilon,i},$$
which, by (51) implies, up to subsequences, $v_{\varepsilon_j} \rightharpoonup v$ weakly in $W^{1,p}(\mathbb{R}^N)$. We now prove that $v \neq 0$. Let us set

$$d_j(y) := \begin{cases} V(x_{\varepsilon_j} + \varepsilon_j y) - \dfrac{g(x_{\varepsilon_j} + \varepsilon_j y, v_{\varepsilon_j}(y))}{v_{\varepsilon_j}^{p-1}(y)} & \text{if } v_{\varepsilon_j}(y) \neq 0 \\ 0 & \text{if } v_{\varepsilon_j}(y) = 0, \end{cases}$$

$$A(y, s, \xi) := j_\xi(x_{\varepsilon_j} + \varepsilon_j y, s, \xi),$$
$$B(y, s, \xi) := d_j(y) s^{p-1},$$
$$C(y, s) := j_s(x_{\varepsilon_j} + \varepsilon_j y, s, Dv_{\varepsilon_j}(y))$$

for every $y \in \mathbb{R}^N$, $s \in \mathbb{R}^+$ and $\xi \in \mathbb{R}^N$. Taking into account the growth of condition on $j_\xi$, the strict convexity of $j$ in $\xi$ and condition (8), we get

$$A(y, s, \xi) \cdot \xi \geqslant \nu |\xi|^p,$$
$$|A(y, s, \xi)| \leqslant c_2 |\xi|^{p-1},$$
$$|B(y, s, \xi)| \leqslant |d_j(y)| |s|^{p-1}.$$

Moreover, by condition (10) we have
$$s \geqslant R \implies C(y, s) s \geqslant 0$$
for every $y \in \mathbb{R}^N$ and $s \in \mathbb{R}^+$. By the growth of conditions on $g$, we have that for $\delta$ sufficiently small $d_j \in L^{\frac{N}{p-\delta}}(B_{2\varrho})$ for every $\varrho > 0$ and
$$S = \sup_j \|d_j\|_{L^{\frac{N}{p-\delta}}(B_{2\varrho})} \leqslant D \sup_{j \in \mathbb{N}} \|v_{\varepsilon_j}\|_{L^{p^*}(B_{2\varrho})} < +\infty$$
for some $D > 0$. Since we have $\mathrm{div}(A(y, v_{\varepsilon_j}, Dv_{\varepsilon_j})) = B(y, v_{\varepsilon_j}, Dv_{\varepsilon_j}) + C(y, v_{\varepsilon_j})$ for every $j \in \mathbb{N}$, by virtue of [26, Theorem 1 and Remark at p.261] there exists a radius $\varrho > 0$ and a positive constant $M = M(\nu, c_2, S\varrho^\delta)$ such that
$$\sup_{j \in \mathbb{N}} \max_{y \in B_\varrho} |v_{\varepsilon_j}(y)| \leqslant M(2\varrho)^{-N/p} \sup_{j \in \mathbb{N}} \|v_{\varepsilon_j}\|_{L^p(B_{2\varrho})} < +\infty$$
so that $(v_{\varepsilon_j})$ is uniformly bounded in $B_\varrho$. Then, by [19, Théorème 1.1, Chapitre IV], up to a subsequence $(v_{\varepsilon_j})$ converges uniformly to $v$ in a small neighbourhood of zero. This yields $v(0) = \lim_j v_{\varepsilon_j}(0) = \lim_j w_{\varepsilon_j}(x_{\varepsilon_j}) \geqslant \ell$.

Up to a rototranslation, it is easily seen that $v$ is a positive solution of
$$-\mathrm{div}(j_\xi(\bar{x}, v, Dv)) + j_s(\bar{x}, v, Dv) + V(\bar{x}) v^{p-1} = \mathbb{1}_{\{x_1 < 0\}} f(v) + \mathbb{1}_{\{x_1 > 0\}} \widetilde{f}(v).$$
By Lemma 3.3 it follows that $v$ is actually a nontrivial solution of
$$-\mathrm{div}(j_\xi(\bar{x}, v, Dv)) + j_s(\bar{x}, v, Dv) + V(\bar{x}) v^{p-1} = f(v).$$

Then, by Theorem 3.2 and (22), we have $I_{\bar{x}}(v) = c_{\bar{x}} \geqslant c_i$. In order to conclude the proof, it is sufficient to prove that

$$(54) \quad \liminf_j \varepsilon_j^{-N} d_{\varepsilon_j,i} = \liminf_j \varepsilon_j^{-N} J_{\varepsilon_j,i}(w_{\varepsilon_j}) \geqslant I_{\bar{x}}(v).$$

STEP II. We prove (54). It results

$$\varepsilon_j^{-N} J_{\varepsilon_j,i}(w_{\varepsilon_j}) = \int_{\hat{\Lambda}_{\varepsilon_j,i}} j(x_{\varepsilon_j} + \varepsilon_j y, v_{\varepsilon_j}, Dv_{\varepsilon_j}) +$$
$$+ \frac{1}{p}\int_{\hat{\Lambda}_{\varepsilon_j,i}} V(x_{\varepsilon_j} + \varepsilon_j y) v_{\varepsilon_j}^p - \int_{\hat{\Lambda}_{\varepsilon_j,i}} G(x_{\varepsilon_j} + \varepsilon_j y, v_{\varepsilon_j})$$

where $\hat{\Lambda}_{\varepsilon_j,i} = \{y \in \mathbb{R}^N : x_{\varepsilon_j} + \varepsilon_j y \in \hat{\Lambda}_i\}$. By Lemma 2.2, we have $v_{\varepsilon_j} \to v$ strongly in $W^{1,p}_{\mathrm{loc}}(\mathbb{R}^N)$. Following the same computations of Theorem 3.2, Step II, we deduce that for all $\delta > 0$ there exists $\bar{\varrho} > 0$ such that for all $\varrho > \bar{\varrho}$ we have

$$\liminf_j \int_{\hat{\Lambda}_{\varepsilon_j,i} \setminus B_\varrho} \Big[ j(x_{\varepsilon_j} + \varepsilon_j y, v_{\varepsilon_j}, Dv_{\varepsilon_j}) +$$
$$+ \frac{1}{p} V(x_{\varepsilon_j} + \varepsilon_j y) v_{\varepsilon_j}^p - G(x_{\varepsilon_j} + \varepsilon_j y, v_{\varepsilon_j}) \Big] \geqslant -\delta.$$

Futhermore we have

$$\lim_j \int_{B_\varrho} \Big[ j(x_{\varepsilon_j} + \varepsilon_j y, v_{\varepsilon_j}, Dv_{\varepsilon_j}) +$$
$$+ \frac{1}{p} V(x_{\varepsilon_j} + \varepsilon_j y) v_{\varepsilon_j}^p - G(x_{\varepsilon_j} + \varepsilon_j y, v_{\varepsilon_j}) \Big] = I_{\bar{x}}(v, B_\varrho),$$

where

$$I_{\bar{x}}(v, B_\varrho) := \int_{B_\varrho} \Big[ j(\bar{x}, v, Dv) + \frac{1}{p} V(\bar{x}) v^p - F(v) \Big].$$

We conclude that for all $\varrho > \bar{\varrho}$

$$\liminf_j \varepsilon_j^{-N} J_{\varepsilon_j,i}(w_{\varepsilon_j}) \geqslant I_{\bar{x}}(v, B_\varrho) - \delta,$$

and (54) follows letting $\varrho \to +\infty$ and $\delta \to 0$. $\square$

Let us now consider the class

$$\Gamma_\varepsilon := \Big\{ \gamma \in C([0,1]^k, W_V(\Omega)) : \ \gamma \text{ satisfies conditions } (a), (b), (c), (d) \Big\},$$

where:

(a) $\gamma(t) = \sum_{i=1}^k \gamma_i(t_i)$ for every $t \in \partial[0,1]^k$, with $\gamma_i \in C([0,1], W_V(\Omega))$;

(b) $\mathrm{supt}(\gamma_i(t_i)) \subset \Lambda_i$ for every $t_i \in [0,1]$ and $i = 1, \ldots, k$;

(c) $\gamma_i(0) = 0$ and $J_\varepsilon(\gamma_i(1)) < 0$ for every $i = 1, \ldots, k$;

(d) $\varepsilon^{-N} E_\varepsilon(\gamma(t)) \leqslant \sum_{i=1}^k c_i + \sigma$ for every $t \in \partial[0,1]^k$,

where $0 < \sigma < \frac{1}{2}\min\{c_i : i = 1,\ldots,k\}$. We set

(55) $$c_\varepsilon := \inf_{\gamma \in \Gamma_\varepsilon} \sup_{t \in [0,1]^k} E_\varepsilon(\gamma(t)).$$

**Lemma 4.2.** *For $\varepsilon$ small enough $\Gamma_\varepsilon \neq \emptyset$ and*

(56) $$\lim_{\varepsilon \to 0^+} \varepsilon^{-N} c_\varepsilon = \sum_{i=1}^{k} c_i.$$

*Proof.* Firstly, let us prove that for $\varepsilon$ small $\Gamma_\varepsilon \neq \emptyset$ and

(57) $$c_\varepsilon \leqslant \varepsilon^N \sum_{i=1}^{k} c_i + o(\varepsilon^N).$$

By definition of $c_i$, for all $\delta > 0$ there exists $\gamma_i \in \mathscr{P}_i$ with

(58) $$c_i \leqslant \max_{\tau \in [0,1]} I_{x_i}(\gamma_i(\tau)) \leqslant c_i + \frac{\delta}{2k}$$

where the $x_i$s are as in (2) and

$$\mathscr{P}_i := \left\{\gamma_i \in C([0,1], W^{1,p}(\mathbb{R}^N)) : \gamma_i(0) = 0, I_{x_i}(\gamma_i(1)) < 0\right\}.$$

We choose $\delta$ so that $\delta < \min\{\sigma, k\sigma_i\}$. Let us set

$$\hat{\gamma}_i(\tau)(x) = \eta_i(x)\gamma_i(\tau)\left(\frac{x - x_i}{\varepsilon}\right) \quad \text{for every } \tau \in [0,1] \text{ and } x \in \Omega,$$

where $\eta_i \in C_c^\infty(\mathbb{R}^N)$, $0 \leqslant \eta_i \leqslant 1$, $\mathrm{supt}(\eta_i) \subseteq \Lambda_i$, and $x_i \in \mathrm{int}(\{\eta_i = 1\})$. We have

(59) $$J_\varepsilon(\hat{\gamma}_i(\tau)) = \int_\Omega \varepsilon^p j(x, \hat{\gamma}_i(\tau), D\hat{\gamma}_i(\tau)) +$$
$$+ \frac{1}{p}\int_\Omega V(x)|\hat{\gamma}_i(\tau)|^p - \int_\Omega G(x, \hat{\gamma}_i(\tau)).$$

Since it results

$$D\hat{\gamma}_i(\tau) = D\eta_i(x)\gamma_i(\tau)\left(\frac{x - x_i}{\varepsilon}\right) + \frac{1}{\varepsilon}\eta_i(x)D\gamma_i(\tau)\left(\frac{x - x_i}{\varepsilon}\right),$$

and for all $\xi_1, \xi_2 \in \mathbb{R}^N$ there exists $t \in [0,1]$ with

$$j(x, s, \xi_1 + \xi_2) = j(x, s, \xi_2) + j_\xi(x, s, t\xi_1 + \xi_2) \cdot \xi_1,$$

taking into account the $p$–homogeneity of $j$, the term

$$\varepsilon^p \int_\Omega j(x, \hat{\gamma}_i(\tau), D\hat{\gamma}_i(\tau))$$

has the same behavior of

(60) $$\int_\Omega j\left(x, \eta_i(x)\gamma_i(\tau)\left(\frac{x - x_i}{\varepsilon}\right), \eta_i(x)D\gamma_i\left(\frac{x - x_i}{\varepsilon}\right)\right)$$

up to an error given by

(61) $$\varepsilon^p \int_\Omega j_\xi(x, s(x), t(x)\xi_1(x) + \xi_2(x)) \cdot \xi_1(x),$$

where we have set

$$s(x) := \hat{\gamma}_i(\tau)(x),$$

$$\xi_1(x) := D\eta_i(x)\gamma_i(\tau)\left(\frac{x - x_i}{\varepsilon}\right),$$

$$\xi_2(x) := \frac{1}{\varepsilon}\eta_i(x)D\gamma_i(\tau)\left(\frac{x - x_i}{\varepsilon}\right),$$

and $t(x)$ is a function with $0 \leqslant t(x) \leqslant 1$ for every $x \in \Omega$.

We proceed in the estimation of (61). We obtain

$$\varepsilon^p \left| \int_\Omega j_\xi(x, s(x), t(x)\xi_1(x) + \xi_2(x)) \cdot \xi_1(x) \right| \leqslant$$

$$\leqslant \widetilde{c}_2 \varepsilon^p \int_\Omega |\xi_1(x)|^p + \widetilde{c}_2 \varepsilon^p \int_\Omega |\xi_2(x)|^{p-1} |\xi_1(x)|.$$

Making the change of variable $y = \frac{x - x_i}{\varepsilon}$, we obtain

$$\varepsilon^p \left| \int_\Omega j_\xi(x, s(x), t(x)\xi_1(x) + \xi_2(x)) \cdot \xi_1(x) \right| \leqslant$$

$$\leqslant \widetilde{c}_2 \varepsilon^{p+N} \int_{\mathbb{R}^N} |D\eta_i(x_i + \varepsilon y)|^p |\gamma_i(\tau)(y)|^p +$$

$$+ \widetilde{c}_2 \varepsilon^{N+1} \int_{\mathbb{R}^N} |\eta_i(x_i + \varepsilon y)|^{p-1} |D\gamma_i(\tau)(y)|^{p-1} \cdot$$

$$\cdot |D\eta_i(x_i + \varepsilon y)| |\gamma_i(\tau)(y)| = o(\varepsilon^N)$$

where $o(\varepsilon^N)$ is independent of $\tau$, since $\gamma_i$ has compact values in $W^{1,p}(\mathbb{R}^N)$. Changing the variable also in (60) yields

$$\int_\Omega j\left(x, \eta_i(x)\gamma_i(\tau)\left(\frac{x - x_i}{\varepsilon}\right), \eta_i(x)D\gamma_i(\tau)\left(\frac{x - x_i}{\varepsilon}\right)\right) =$$

$$= \varepsilon^N \int_{\mathbb{R}^N} j(x_i + \varepsilon y, \eta_i(x_i + \varepsilon y)\gamma_i(\tau)(y), \eta_i(x_i + \varepsilon y)D\gamma_i(\tau)(y)).$$

By the Dominated Convergence Theorem we get

$$\lim_{\varepsilon \to 0} \int_{\mathbb{R}^N} j(x_i + \varepsilon y, \eta_i(x_i + \varepsilon y)\gamma_i(\tau)(y), \eta_i(x_i + \varepsilon y)D\gamma_i(\tau)(y)) =$$

$$= \int_{\mathbb{R}^N} j(x_i, \gamma_i(\tau)(y), D\gamma_i(\tau)(y))$$

uniformly with respect to $\tau$. Reasoning in a similar fashion for the other terms in (59), we conclude that for $\varepsilon$ small enough

(62) $$J_\varepsilon(\hat{\gamma}_i(\tau)) = \varepsilon^N I_{x_i}(\gamma_i(\tau)) + o(\varepsilon^N)$$

for every $\tau \in [0,1]$ with $o(\varepsilon^N)$ independent of $\tau$. Let us now set

$$\gamma_0(\tau_1, \ldots, \tau_k) := \sum_{i=1}^{k} \hat{\gamma}_i(\tau_i).$$

Since $\mathrm{supt}(\hat{\gamma}_i(\tau)) \subseteq \Lambda_i$ for every $\tau$, we have that $J_{\varepsilon,i}(\hat{\gamma}_i(\tau)) = J_\varepsilon(\hat{\gamma}_i(\tau))$; then, by the choice of $\delta$, we get for $\varepsilon$ small

$$[J_{\varepsilon,i}(\hat{\gamma}_i(\tau))_+]^{\frac{1}{2}} - \varepsilon^{\frac{N}{2}}(c_i + \sigma_i)^{\frac{1}{2}} =$$
$$= [J_\varepsilon(\hat{\gamma}_i(\tau))_+]^{\frac{1}{2}} - \varepsilon^{\frac{N}{2}}(c_i + \sigma_i)^{\frac{1}{2}} =$$
$$= \varepsilon^{\frac{N}{2}}[I_{x_i}(\gamma_i(\tau)) + o(1)]^{\frac{1}{2}} - \varepsilon^{\frac{N}{2}}(c_i + \sigma_i)^{\frac{1}{2}} \leqslant$$
$$\leqslant \varepsilon^{\frac{N}{2}}\left[c_i + \frac{\delta}{2k} + o(1)\right]^{\frac{1}{2}} - \varepsilon^{\frac{N}{2}}(c_i + \sigma_i)^{\frac{1}{2}} \leqslant 0,$$

and

$$E_\varepsilon(\gamma_0(\tau_1, \ldots, \tau_k)) = J_\varepsilon(\gamma_0(\tau_1, \ldots, \tau_k)) = \sum_{i=1}^{k} J_\varepsilon(\hat{\gamma}_i(\tau_i)).$$

By (58) and (62) we obtain that for $\varepsilon$ small enough

$$E_\varepsilon(\gamma_0(\tau)) \leqslant \varepsilon^N \sum_{i=1}^{k}\left(c_i + \frac{\delta}{2k}\right) \leqslant \varepsilon^N\left(\sum_{i=1}^{k} c_i + \sigma\right)$$

so that the class $\Gamma_\varepsilon$ is not empty. Moreover, we have

$$\limsup_{\varepsilon \to 0^+} \frac{c_\varepsilon}{\varepsilon^N} \leqslant \sum_{i=1}^{k} c_i + \delta$$

and, by the arbitrariness of $\delta$, we have conclude that (57) holds. Let us now prove that

(63) $$c_\varepsilon \geqslant \varepsilon^N \sum_{i=1}^{k} c_i + o(\varepsilon^N).$$

Given $\gamma \in \Gamma_\varepsilon$, by a variant of [7, Proposition 3.4] there exists $\bar{t} \in [0,1]^k$ such that

$$J_{\varepsilon,i}(\gamma(\bar{t})) \geqslant d_{\varepsilon,i}$$

for all $i = 1, \ldots, k$, where the $d_{\varepsilon,i}$s are as in (50). Then we have by Lemma 4.1

$$\sup_{t \in [0,1]^k} J_\varepsilon(\gamma(t)) \geqslant \sup_{t \in [0,1]^k} \sum_{i=1}^{k} J_{\varepsilon,i}(\gamma(t)) \geqslant \sum_{i=1}^{k} d_{\varepsilon,i} = \varepsilon^N \sum_{i=1}^{k} c_i + o(\varepsilon^N),$$

which implies the assertion. $\square$

**Corollary 4.3.** *For every $\varepsilon > 0$ there exists a critical point $u_\varepsilon \in W_V(\Omega)$ of the functional $E_\varepsilon$ such that $c_\varepsilon = E_\varepsilon(u_\varepsilon)$. Moreover $\|u_\varepsilon\|_{W_V} \to 0$ as $\varepsilon \to 0$.*

*Proof.* By combining Lemma 2.3 with (*b*) of Proposition 6.8 it results that $E_\varepsilon$ satisfies the Palais–Smale condition for every $c \in \mathbb{R}$ (see Definition 6.3). Then, taking into account Lemma 4.2, for every $\varepsilon > 0$ the (nonsmooth) Mountain–Pass Theorem (see [6]) for the class $\Gamma_\varepsilon$ provides the desired critical point $u_\varepsilon$ of $E_\varepsilon$. To prove the second assertion we may argue as in Step I of Lemma 2.3 with $u_h$ replaced by $u_\varepsilon$ and $c$ replaced by $E_\varepsilon(u_\varepsilon)$. Thus, from inequality (32), for every $\varepsilon > 0$ we get

$$(64) \quad \int_\Omega \left(|Du_\varepsilon|^p + V(x)|u_\varepsilon|^p\right) \leqslant \frac{\vartheta E_\varepsilon(u_\varepsilon) + d\varepsilon^N}{\min\left\{(\vartheta - \gamma - p)\nu\varepsilon^p, \frac{\vartheta}{p} - \frac{\vartheta}{2\kappa} - 1\right\}}.$$

By virtue of Lemma 4.2, this yields

$$\int_\Omega \left(|Du_\varepsilon|^p + V(x)|u_\varepsilon|^p\right) \leqslant \left\{\frac{\vartheta(c_1 + \cdots + c_k) + d}{(\vartheta - \gamma - p)\nu}\right\}\varepsilon^{N-p} + o(\varepsilon^{N-p}),$$

as $\varepsilon \to 0$, which implies the assertion. □

Let us now set:

$$\Omega_\varepsilon := \{y \in \mathbb{R}^N : \varepsilon y \in \Omega\}, \qquad v_\varepsilon(y) := u_\varepsilon(\varepsilon y) \in W^{1,p}(\Omega_\varepsilon),$$
$$\hat{\Lambda}_{\varepsilon,i} := \{y \in \mathbb{R}^N : \varepsilon y \in \hat{\Lambda}_i\}, \qquad \Lambda_\varepsilon := \{y \in \mathbb{R}^N : \varepsilon y \in \Lambda\}.$$

**Lemma 4.4.** *The function $v_\varepsilon$ is a solution of the equation*

$$(65) \quad -\operatorname{div}\left((1+\theta_\varepsilon(\varepsilon y))j_\xi(\varepsilon y, v, Dv)\right) + (1+\theta_\varepsilon(\varepsilon y))j_s(\varepsilon y, v, Dv) +$$
$$+ (1+\theta_\varepsilon(\varepsilon y))V(\varepsilon y)v^{p-1} = (1+\theta_\varepsilon(\varepsilon y))g(\varepsilon y, v) \quad \text{in } \Omega_\varepsilon,$$

*where for every $\varepsilon > 0$*

$$(66) \quad \theta_\varepsilon(x) := \sum_{i=1}^k \theta_{\varepsilon,i}\mathbb{1}_{\hat{\Lambda}_i}(x), \quad \theta_{\varepsilon,i} \in [0, M],$$

$$\theta_{\varepsilon,i} := M\left[(J_{\varepsilon,i}(u_\varepsilon)_+)^{1/2} - \varepsilon^{N/2}(c_i + \sigma_i)^{1/2}\right]_+ (J_{\varepsilon,i}(u_\varepsilon)_+)^{-1/2}.$$

*Proof.* It suffices to expand $E'_\varepsilon(u_\varepsilon)(\varphi) = 0$ for every $\varphi \in C^\infty_c(\Omega)$. □

**Corollary 4.5.** *The sequence $(v_\varepsilon)$ is bounded in $W^{1,p}(\mathbb{R}^N)$.*

*Proof.* It suffices to combine Lemma 4.2 with the inequality

$$\int_{\mathbb{R}^N} \left(|Dv_\varepsilon|^p + V(x)|v_\varepsilon|^p\right) \leqslant \frac{\vartheta\varepsilon^{-N}c_\varepsilon + d}{\min\left\{(\vartheta - \gamma - p)\nu, \frac{\vartheta}{p} - \frac{\vartheta}{2\kappa} - 1\right\}}$$

which follows by (64). □

The following lemma "kills" the second penalization term of $E_\varepsilon$.

**Lemma 4.6.** *We have*

$$(67) \qquad \lim_{\varepsilon \to 0} \varepsilon^{-N} J_{\varepsilon,i}(u_\varepsilon) = c_i$$

*for every $i = 1, \ldots, k$.*

*Proof.* Let us first prove that, as $\varrho \to +\infty$,

$$\text{(68)} \quad \limsup_{\varepsilon \to 0^+} \int_{\Omega_\varepsilon \setminus \mathcal{N}_\varrho(\Lambda_\varepsilon)} \left(|Dv_\varepsilon|^p + |v_\varepsilon|^p\right) = o(1),$$

where $\mathcal{N}_\varrho(\Lambda_\varepsilon) := \{y \in \mathbb{R}^N : \text{dist}(y, \Lambda_\varepsilon) < \varrho\}$. By Proposition 6.4, we can test equation (65) with $\psi_{\varepsilon,\varrho} v_\varepsilon e^{\zeta(v_\varepsilon)}$, where $\psi_{\varepsilon,\varrho} = 1 - \sum_{i=1}^k \psi^i_{\varepsilon,\varrho}$, $\psi^i_{\varepsilon,\varrho} \in C^\infty(\mathbb{R}^N)$,

$$\psi^i_{\varepsilon,\varrho} = 1 \text{ if } \text{dist}(y, \Lambda_{\varepsilon,i}) < \varrho/2,$$

$$\psi^i_{\varepsilon,\varrho} = 0 \text{ if } \text{dist}(y, \Lambda_{\varepsilon,i}) > \varrho$$

and the function $\zeta$ is defined as in (35). By virtue of (1), (7), the boundedness of $(v_\varepsilon)$ in $W^{1,p}(\mathbb{R}^N)$ and (37) there exist $C, C' > 0$ such that

$$C \int_{\Omega_\varepsilon \setminus \mathcal{N}_\varrho(\Lambda_\varepsilon)} \left(|Dv_\varepsilon|^p + |v_\varepsilon|^p\right) \leqslant$$

$$\leqslant \int_{\Omega_\varepsilon \setminus \Lambda_\varepsilon} (1 + \theta_\varepsilon(\varepsilon y)) \left[pj(\varepsilon y, v_\varepsilon, Dv_\varepsilon) + \left\{V(\varepsilon y) - \frac{\widetilde{f}(v_\varepsilon)}{v_\varepsilon^{p-1}}\right\} v_\varepsilon^p\right] \psi_{\varepsilon,\varrho} e^{\zeta(v_\varepsilon)} =$$

$$= -\int_{\Omega_\varepsilon \setminus \Lambda_\varepsilon} (1 + \theta_\varepsilon(\varepsilon y))[j_s(\varepsilon y, v_\varepsilon, Dv_\varepsilon) + p\zeta'(v_\varepsilon)j(\varepsilon y, v_\varepsilon, Dv_\varepsilon)]v_\varepsilon \psi_{\varepsilon,\varrho} e^{\zeta(v_\varepsilon)} +$$

$$- \int_{\Omega_\varepsilon \setminus \Lambda_\varepsilon} (1 + \theta_\varepsilon(\varepsilon y))j_\xi(\varepsilon y, v_\varepsilon, Dv_\varepsilon) \cdot D\psi_{\varepsilon,\varrho} v_\varepsilon e^{\zeta(v_\varepsilon)} \leqslant$$

$$\leqslant 2e^{\bar{M}R} \int_{\Omega_\varepsilon \setminus \Lambda_\varepsilon} |D\psi_{\varepsilon,\varrho}||j_\xi(\varepsilon y, v_\varepsilon, Dv_\varepsilon)|v_\varepsilon \leqslant \frac{\widetilde{C}}{\varrho}\|Dv_\varepsilon\|_p^{p-1}\|v_\varepsilon\|_p \leqslant \frac{C'}{\varrho},$$

which implies (68). Now, to prove (67), we adapt the argument of [12, Lemma 2.1] to our context. It is sufficient to prove that

$$\text{(69)} \quad \lim_{\varepsilon \to 0} \varepsilon^{-N} J_{\varepsilon,i}(u_\varepsilon) \leqslant c_i + \sigma_i$$

for every $i = 1, \ldots, k$. Then (67) follows by arguing exactly as in [12, Lemma 2.4]. By contradiction, let us suppose that for some $\varepsilon_j \to 0$ we have

$$\text{(70)} \quad \limsup_j \varepsilon_j^{-N} J_{\varepsilon_j,i}(u_{\varepsilon_j}) > c_i + \sigma_i.$$

Then there exists $\lambda > 0$ with

$$\int_{\hat{\Lambda}_{\varepsilon_j,i}} \left(|Dv_{\varepsilon_j}|^p + |v_{\varepsilon_j}|^p\right) \geqslant \lambda,$$

and so by (68) there exists $\varrho > 0$ such that for $j$ large enough

$$\int_{\mathcal{N}_\varrho(\Lambda_{\varepsilon_j,i})} \left(|Dv_{\varepsilon_j}|^p + |v_{\varepsilon_j}|^p\right) \geqslant \frac{\lambda}{2}.$$

Following [12, Lemma 2.1], P.L. Lions' concentration compactness argument [21] yields the existence of $S > 0$, $\rho > 0$ and a sequence $y_j \in \Lambda_{\varepsilon_j, i}$ such that for $j$ large enough

$$\int_{B_S(y_j)} v_{\varepsilon_j}^p \geqslant \rho. \tag{71}$$

Let us set $v_j(y) := v_{\varepsilon_j}(y_j + y)$, and let $\varepsilon_j y_j \to \bar{x} \in \Lambda_i$. By Corollary 4.5, we may assume that $v_j$ weakly converges to some $v$ in $W^{1,p}(\mathbb{R}^N)$. By Lemma 2.2, we have that $v_j \to v$ strongly in $W^{1,p}_{\text{loc}}(\mathbb{R}^N)$; note that $v \neq 0$ by (71). In the case $\text{dist}(y_j, \partial \Lambda_{\varepsilon_j, i}) \to +\infty$, since $v_j$ satisfies in $-y_j + \Lambda_{\varepsilon_j, i}$ the equation

$$-\text{div}(j_\xi(\varepsilon_j y_j + \varepsilon_j y, v_j, Dv_j)) + j_s(\varepsilon_j y_j + \varepsilon_j y, v_j, Dv_j) + \\ + V(\varepsilon_j y_j + \varepsilon_j y) v^{p-1} = f(v_j),$$

$v$ satisfies on $\mathbb{R}^N$ the equation

$$-\text{div}(j_\xi(\bar{x}, v, Dv)) + j_s(\bar{x}, v, Dv) + V(\bar{x}) v^{p-1} = f(v). \tag{72}$$

If $\text{dist}(y_j, \partial \Lambda_{\varepsilon_j, i}) \leqslant C < +\infty$, we deduce that $v$ satisfies an equation of the form (45), and by Lemma 3.3, we conclude that $v$ satisfies equation (72). Since $v$ is a nontrivial critical point for $I_{\bar{x}}$, by (11) and Theorem 3.2, recalling that $c_i \leqslant c_{\bar{x}} \leqslant c_i + \sigma_i$, we get $c_i \leqslant I_{\bar{x}}(v) \leqslant c_i + \sigma_i$. Then we can find a sequence $R_j \to +\infty$ such that

$$\lim_j \int_{B_{R_j}(y_j)} j(\varepsilon_j y, v_{\varepsilon_j}, Dv_{\varepsilon_j}) + \frac{1}{p} V(\varepsilon_j y) |v_{\varepsilon_j}|^p - G(\varepsilon_j y, v_{\varepsilon_j}) = \\ = I_{\bar{x}}(v) \leqslant c_i + \sigma_i.$$

Then by (70) we deduce that for $j$ large enough

$$\int_{\hat{\Lambda}_{\varepsilon_j, i} \setminus B_{R_j}(y_j)} \left( |Dv_{\varepsilon_j}|^p + |v_{\varepsilon_j}|^p \right) \geqslant \lambda > 0.$$

Reasoning as above, there exist $\widetilde{S}, \widetilde{\rho} > 0$ and a sequence $\widetilde{y}_j \in \Lambda_{\varepsilon_j, i} \setminus B_{R_j}(y_j)$ such that

$$\int_{B_{\widetilde{S}}(\widetilde{y}_j)} v_{\varepsilon_j}^p \geqslant \widetilde{\rho} > 0. \tag{73}$$

Let $\varepsilon_j \widetilde{y}_j \to \widetilde{x} \in \Lambda_i$; then we have $\widetilde{v}_j(y) := v_{\varepsilon_j}(\widetilde{y}_j + y) \to \widetilde{v}$ weakly in $W^{1,p}(\mathbb{R}^N)$, where $\widetilde{v}$ is a nontrivial solution of the equation

$$-\text{div}(j_\xi(\widetilde{x}, v, Dv)) + j_s(\widetilde{x}, v, Dv) + V(\widetilde{x}) v^{p-1} = f(v).$$

As before we get $I_{\widetilde{x}}(\widetilde{v}) \geqslant c_i$. We are now in a position to deduce that

$$\liminf_j \varepsilon_j^{-N} J_{\varepsilon_j, i}(u_\varepsilon) > I_{\bar{x}}(v) + I_{\widetilde{x}}(\widetilde{v}) \geqslant 2c_i.$$

In fact, $v_{\varepsilon_j}$ satisfies in $\hat{\Lambda}_{\varepsilon_j, i}$ the equation

$$-\text{div}(j_\xi(\varepsilon_j y, v_{\varepsilon_j}, Dv_{\varepsilon_j})) + j_s(\varepsilon_j y, v_{\varepsilon_j}, Dv_{\varepsilon_j}) + V(\varepsilon_j y) v_{\varepsilon_j}^{p-1} = g(\varepsilon_j y, v_{\varepsilon_j}). \tag{74}$$

Since $y_j, \widetilde{y}_j \in \Lambda_{\varepsilon_j, i}$, for $j$ large enough $B_{j,R} := B(y_j, R) \cup B(\widetilde{y}_j, R) \subset \hat{\Lambda}_{\varepsilon_j, i}$, and so we can test (74) with

$$\varphi(y) = \left[\psi\left(\frac{|y - y_j|}{R}\right) + \psi\left(\frac{|y - \widetilde{y}_j|}{R}\right) - 1\right] v_{\varepsilon_j}(y)$$

where $\psi \in C^\infty(\mathbb{R})$ with $0 \leqslant \psi \leqslant 1$, $\psi(s) = 0$ for $s \leqslant 1$ and $\psi(s) = 1$ for $s \geqslant 2$. Reasoning as in Lemma 4.1, we have that for all $\delta > 0$ there exists $\bar{R}$ such that for all $R > \bar{R}$ we have

$$\int_{\hat{\Lambda}_{\varepsilon_j,i} \setminus B_{j,R}} \left[j(\varepsilon_j y, v_{\varepsilon_j}, Dv_{\varepsilon_j}) + \frac{1}{p} V(\varepsilon_j y) |v_{\varepsilon_j}|^p - G(\varepsilon_j y, v_{\varepsilon_j})\right] \geqslant -\delta$$

so that

$$\liminf_j \varepsilon_j^{-N} J_{\varepsilon_j, i}(u_{\varepsilon_j}) \geqslant I_{\bar{x}}(v, B_R) + I_{\widetilde{x}}(\widetilde{v}, B_R) - \delta.$$

Letting $R \to +\infty$ and $\delta \to 0$, we get

(75) $$\liminf_j \varepsilon_j^{-N} J_{\varepsilon_j, i}(u_{\varepsilon_j}) > 2c_i.$$

The same arguments apply to the functional $J_\varepsilon$: we have that

(76) $$\liminf_j \varepsilon_j^{-N} J_{\varepsilon_j}(u_{\varepsilon_j}) \geqslant 2c_i.$$

Then by combining (75) and (76) we obtain

$$\liminf_j \varepsilon_j^{-N} E_{\varepsilon_j}(u_{\varepsilon_j}) \geqslant 2c_i + M\left[(2c_i)^{\frac{1}{2}} - (c_i + \sigma_i)^{\frac{1}{2}}\right]_+^2.$$

By Lemma 4.2, we have

$$M\left[(2c_i)^{\frac{1}{2}} - (c_i + \sigma_i)^{\frac{1}{2}}\right]_+^2 \leqslant \sum_{i=1}^k c_i,$$

against the choice of $M$. $\square$

## 5. Proofs of the main results

We are now ready to prove the main results of the paper.

PROOF OF THEOREM 1.2. Let us consider the sequence $(u_\varepsilon)$ of critical points of $E_\varepsilon$ given by Corollary 4.3. We have that $\|u_\varepsilon\|_{W_V} \to 0$. Since $u_\varepsilon$ satisfies

$$-\operatorname{div}\left((1 + \theta_\varepsilon(x)) j_\xi(x, v, Dv)\right) + (1 + \theta_\varepsilon(x)) j_s(x, v, Dv) +$$
$$+ (1 + \theta_\varepsilon(x)) V(x) v^{p-1} = (1 + \theta_\varepsilon(x)) g(x, v) \quad \text{in } \Omega,$$

with $\theta_\varepsilon$ defined as in (66), by the regularity results of [19] $u_\varepsilon$ is locally Hölder continuous in $\Omega$. We claim that there exists $\sigma > 0$ such that

(77) $$u_\varepsilon(x_{\varepsilon, i}) = \sup_{\Lambda_i} u_\varepsilon > \sigma > 0$$

for every $\varepsilon$ sufficiently small and $i = 1, \ldots, k$: moreover

(78) $$\lim_{\varepsilon \to 0} \operatorname{dist}(x_{\varepsilon, i}, \mathcal{M}_i) = 0$$

for every $i = 1, \ldots, k$, where the $\mathscr{M}_i$s are the sets of minima of $V$ in $\Lambda_i$. In fact, let us assume that there exists $i_0 \in \{1, \ldots, k\}$ such that $u_\varepsilon(x_{\varepsilon,i_0}) \to 0$ as $\varepsilon \to 0$. Therefore, $u_\varepsilon \to 0$ uniformly on $\Lambda_{i_0}$ as $\varepsilon \to 0$, which implies that

(79) $$\sup_{y \in \Lambda_{\varepsilon,i_0}} v_\varepsilon(y) \to 0 \quad \text{as } \varepsilon \to 0,$$

where $v_\varepsilon(y) := u_\varepsilon(\varepsilon y)$. On the other hand, since by (67) we have

$$\lim_{\varepsilon \to 0} \varepsilon^{-N} J_{\varepsilon, i_0}(u_\varepsilon) = c_{i_0} > 0,$$

considering $\widetilde{\Lambda}_{i_0}$ relatively compact in $\Lambda_{i_0}$, following the proof of Lemma 4.6, we find $S > 0$ and $\varrho > 0$ such that

$$\sup_{y \in \widetilde{\Lambda}_{\varepsilon,i_0}} \int_{B_S(y)} v_\varepsilon^p \geqslant \varrho$$

for every $\varepsilon \in (0, \varepsilon_0)$, which contradicts (79). We conclude that (77) holds. In order to prove (78), it is sufficient to prove that

$$\lim_{\varepsilon \to 0} V(x_{\varepsilon,i}) = \min_{\Lambda_i} V$$

for every $i = 1, \ldots, k$. Assume by contradiction that for some $i_0$

$$\lim_{\varepsilon \to 0} V(x_{\varepsilon,i_0}) > \min_{\Lambda_{i_0}} V = b_{i_0}.$$

Then, up to a subsequence, $x_{\varepsilon_j,i_0} \to x_{i_0} \in \Lambda_{i_0}$ and $V(x_{i_0}) > b_{i_0}$. Then, arguing as in the proof of Lemma 4.6 and using Theorem 3.2, we would get

$$\liminf_j \varepsilon_j^{-N} J_{\varepsilon_j,i_0}(u_{\varepsilon_j}) \geqslant I_{x_{i_0}}(v) = c_{x_{i_0}} > c_{i_0}$$

which is impossible, in view of (67).

We now prove that

(80) $$\lim_{\varepsilon \to 0} u_\varepsilon = 0 \quad \text{uniformly on } \Omega \setminus \bigcup_{i=1}^k \operatorname{int}(\Lambda_i).$$

Let us first prove that

$$\limsup_{\varepsilon \to 0} \sup_{\partial \Lambda_i} u_\varepsilon = 0 \quad \text{for every } i = 1, \ldots, k.$$

By contradiction, let $i_0 \in \{1, \ldots, k\}$ and $\sigma > 0$ with $u_{\varepsilon_j}(x_j) \geqslant \sigma$ for $(x_j) \subset \partial \Lambda_{i_0}$. Up to a subsequence, $x_j \to x_0 \in \partial \Lambda_{i_0}$. Therefore, taking into account Lemma 3.3 and the local regularity estimates of [26] (see also the end of Step I of Lemma 4.1), the sequence $v_j(y) := u_{\varepsilon_j}(x_j + \varepsilon_j y)$ converges weakly to a nontrivial solution $v \in W^{1,p}(\mathbb{R}^N)$ of

$$-\operatorname{div}(j_\xi(x_0, v, Dv)) + j_s(x_0, v, Dv) + V(x_0)v^{p-1} = f(v) \quad \text{in } \mathbb{R}^N.$$

As $V(x_0) > V(x_{i_0})$, we have

$$\liminf_j \varepsilon_j^{-N} J_{\varepsilon_j, i_0}(u_{\varepsilon_j}) \geqslant I_{x_0}(v) > c_{i_0},$$

which violates (67). Testing the equation with

$$(u_\varepsilon - \max_i \sup_{\partial \Lambda_i} u_\varepsilon)^+ \mathbb{1}_{\Omega \setminus \Lambda} e^{\zeta(u_\varepsilon)},$$

as in Lemma 3.3, this yields that $u_\varepsilon(x) \leqslant \max_i \sup_{\partial \Lambda_i} u_\varepsilon$ for every $x \in \Omega \setminus \Lambda$, so that (80) holds.

By Proposition 2.1, $u_\varepsilon$ is actually a solution of the original problem $(P_\varepsilon)$ because the penalization terms are neutralized by the facts $J_{\varepsilon,i}(u_\varepsilon) < c_i + \sigma_i$ and $u_\varepsilon < \ell$ on $\Omega \setminus \Lambda$ for $\varepsilon$ small. By regularity results, it follows $u_\varepsilon \in C^{1,\beta}_{\text{loc}}(\Omega)$, and so point $(a)$ is proved. Taking into account (77) and (80), we get that $u_\varepsilon$ has a maximum $\bar{x}_\varepsilon \in \Omega$ which coincides with one of the $x_{\varepsilon,i}$s. Considering $\bar{v}_\varepsilon(y) := u_\varepsilon(x_{\varepsilon,i} + \varepsilon y)$, since $\bar{v}_\varepsilon$ is uniformly bounded in $W^{1,p}_{\text{loc}}(\mathbb{R}^N)$, by the local regularity estimates [26], there exists $\sigma'$ with

$$u_\varepsilon(x_{\varepsilon,i}) \leqslant \sigma'$$

for all $i = 1, \ldots, k$. In view of (77), (78) and Corollary 4.3, we conclude that points $(b)$ and $(d)$ are proved. Let us now come to point $(c)$. Let us assume by contradiction that there exists $\bar{r}$, $\delta$, $i_0$ and $\varepsilon_j \to 0$ such that there exists $y_j \in \Lambda_{i_0} \setminus B_{\bar{r}}(x_{\varepsilon_j,i_0})$ with

$$\limsup_j u_{\varepsilon_j}(y_j) \geqslant \delta.$$

We may assume that $y_j \to \bar{y}$, $x_{\varepsilon_j,i_0} \to \bar{x}$, and $\bar{v}_j(y) := u_{\varepsilon_j}(y_j + \varepsilon_j y) \to \bar{v}$, $v_j(y) := u_{\varepsilon_j}(x_{\varepsilon_j,i_0} + \varepsilon_j y) \to v$ strongly in $W^{1,p}_{\text{loc}}(\mathbb{R}^N)$: then, arguing as in Lemma 4.6, it turns out that

$$\liminf_j \varepsilon_j^{-N} J_{\varepsilon_j, i_0}(u_{\varepsilon_j}) \geqslant I_{\bar{x}}(v) + I_{\bar{y}}(\bar{v}) \geqslant 2c_{i_0}$$

which is against (67). We conclude that point $(c)$ holds, and the proof is concluded. □

PROOF OF THEOREM 1.1. If $1 < p \leqslant 2$ and $p < q < p^*$, the equation

(81) $$-\Delta_p u + V(\bar{x}) u^{p-1} = u^{q-1} \quad \text{in } \mathbb{R}^N$$

admits a unique positive $C^1$ solution (up to translations).

Indeed, a solution $u \in C^1(\mathbb{R}^N)$ of (81) exists by Theorem 3.2. By [20, Theorem 1] we have $u(x) \to 0$ as $|x| \to \infty$. Moreover, by [8, Theorem 1.1], the solution $u$ is radially symmetric about some point $x_0 \in \mathbb{R}^N$ and radially decreasing. Then $u$ is a radial ground state solution of (81). By [27, Theorem 1], $u$ is unique (up to translations). Then (11) is satisfied and the assertions follow by Theorem 1.2 applied to the functions $j(x, s, \xi) = \frac{1}{p}|\xi|^p$ and $f(s) = s^{q-1}$. □

6. Appendix. Recalls of nonsmooth critical point theory

In this section we quote from [4, 6] some tools of nonsmooth critical point theory which we use in the paper.

Let us first recall the definition of weak slope for a continuous function.

**Definition 6.1.** Let $X$ be a complete metric space, $f : X \to \mathbb{R}$ be a continuous function, and $u \in X$. We denote by $|df|(u)$ the supremum of the real numbers $\sigma \geqslant 0$ such that there exist $\delta > 0$ and a continuous map
$$\mathscr{H} : B(u,\delta) \times [0,\delta] \to X,$$
such that, for every $v$ in $B(u,\delta)$, and for every $t$ in $[0,\delta]$ it results
$$d(\mathscr{H}(v,t),v) \leqslant t, \quad f(\mathscr{H}(v,t)) \leqslant f(v) - \sigma t.$$
The extended real number $|df|(u)$ is called the weak slope of $f$ at $u$.

The previous notion allows us to give the following definitions.

**Definition 6.2.** We say that $u \in X$ is a critical point of $f$ if $|df|(u) = 0$. We say that $c \in \mathbb{R}$ is a critical value of $f$ if there exists a critical point $u \in X$ of $f$ with $f(u) = c$.

**Definition 6.3.** Let $c \in \mathbb{R}$. We say that $f$ satisfies the Palais–Smale condition at level $c$ ($(PS)_c$ in short), if every sequence $(u_h)$ in $X$ such that $|df|(u_h) \to 0$ and $f(u_h) \to c$ admits a subsequence converging in $X$.

Let us now return to the concrete setting and choose $X = W_V(\Omega)$. Let $\varepsilon > 0$ and consider the functional $f : W_V(\Omega) \to \mathbb{R}$ defined by setting

$$(82) \qquad f(u) = \varepsilon^p \int_\Omega j(x,u,Du) + \frac{1}{p} \int_\Omega V(x)|u|^p - \int_\Omega G(x,u)$$

where $G(x,s) = \int_0^s g(x,t)\,dt$ and $g : \Omega \times \mathbb{R} \to \mathbb{R}$ is now any Carathéodory function. Although $f$ is a nonsmooth functional, its directional derivatives exist along some special directions.

**Proposition 6.4.** Let $u, \varphi \in W_V(\Omega)$ be such that $[j_s(x,u,Du)\varphi]^- \in L^1(\Omega)$. Then we have $j_s(x,u,Du)\varphi \in L^1(\Omega)$, the directional derivative $f'(u)(\varphi)$ exists and
$$f'(u)(\varphi) = \varepsilon^p \int_\Omega j_\xi(x,u,Du) \cdot D\varphi + \varepsilon^p \int_\Omega j_s(x,u,Du)\varphi +$$
$$+ \int_\Omega V(x)|u|^{p-2}u\varphi - \int_\Omega g(x,u)\varphi.$$
In particular, if condition (10) holds, for every $\varphi \in L^\infty(\Omega)$, $\varphi \geqslant 0$ we have $j_s(x,u,Du)\varphi u \in L^1(\Omega)$ and the derivative $f'(u)(\varphi u)$ exists.

**Definition 6.5.** We say that $u$ is a (weak) solution of the problem
$$(83) \quad \begin{cases} -\varepsilon^p \mathrm{div}(j_\xi(x,u,Du)) + \varepsilon^p j_s(x,u,Du) + V(x)|u|^{p-2}u = g(x,u) & \text{in } \Omega \\ u = 0 & \text{on } \partial\Omega \end{cases}$$

if $u \in W_V(\Omega)$ and
$$-\varepsilon^p \mathrm{div}(j_\xi(x,u,Du)) + \varepsilon^p j_s(x,u,Du) + V(x)|u|^{p-2}u = g(x,u)$$
is satisfied in $\mathscr{D}'(\Omega)$.

We now introduce a useful variant of the classical Palais–Smale condition.

**Definition 6.6.** Let $\varepsilon > 0$ and $c \in \mathbb{R}$. We say that $(u_h) \subset W_V(\Omega)$ is a concrete Palais–Smale sequence at level $c$ ($(CPS)_c$–sequence, in short) for $f$, if $f(u_h) \to c$ and
$$j_s(x, u_h, Du_h) \in (W_V(\Omega))' \quad \text{as } h \to +\infty,$$
$$-\varepsilon^p \mathrm{div}(j_\xi(x,u_h,Du_h)) + \varepsilon^p j_s(x,u_h,Du_h) + V(x)|u_h|^{p-2}u_h - g(x,u_h) \to 0$$
strongly in $(W_V(\Omega))'$. We say that $f$ satisfies the concrete Palais–Smale condition at level $c$ ($(CPS)_c$ condition), if every $(CPS)_c$–sequence for $f$ admits a strongly convergent subsequence in $W_V(\Omega)$.

**Proposition 6.7.** *Let $\varepsilon > 0$. Then for every $u$ in $W_V(\Omega)$ with $|df|(u) < +\infty$ we have*
$$-\varepsilon^p \mathrm{div}(j_\xi(x,u,Du)) + \varepsilon^p j_s(x,u,Du) \in (W_V(\Omega))'$$
*and setting*
$$w_u^\varepsilon = -\varepsilon^p \mathrm{div}(j_\xi(x,u,Du)) + \varepsilon^p j_s(x,u,Du) + V(x)|u|^{p-2}u - g(x,u)$$
*it results $\|w_u^\varepsilon\|_{(W_V(\Omega))'} \leqslant |dJ_\varepsilon|(u)$.*

As a consequence of the previous proposition we have the following result.

**Proposition 6.8.** *Let $u \in W_V(\Omega)$, $c \in \mathbb{R}$ and let $(u_h) \subset W_V(\Omega)$.*
*Then the following facts hold:*
  *(a) if $u$ is a critical point of $f$, then $u$ is a weak solution of (83);*
  *(b) if $f$ satisfies the $(CPS)_c$ condition, then $f$ satisfies the $(PS)_c$ condition.*

For suitable versions of the Mountain–Pass Theorem in the nonsmooth framework we refer the reader to [4].

Scuola Internazionale Superiore di Studi Avanzati
Via Beirut 4, I–34014 Trieste, Italy.
 *E-mail address*: `giacomin@sissa.it`

Università Cattolica del Sacro Cuore
Via Musei 41, I–25121 Brescia, Italy.
 *E-mail address*: `m.squassina@dmf.unicatt.it`
 *URL*: `www.dmf.unicatt.it/squassina`